



\magnification 1200
\hsize 13.2cm
\vsize 20cm
\parskip 3pt plus 1pt
\parindent 5mm

\def\\{\hfil\break}


\font\seventeenbf=cmbx10 at 17.28pt

\font\twelvebf=cmbx10 at 12pt
\font\eightbf=cmbx8
\font\sixbf=cmbx6

\font\eighti=cmmi8
\font\sixi=cmmi6

\font\eightrm=cmr8
\font\sixrm=cmr6

\font\eightsy=cmsy8
\font\sixsy=cmsy6

\font\eightit=cmti8
\font\eighttt=cmtt8
\font\eightsl=cmsl8

\font\seventeenbsy=cmbsy10 at 17.28pt

\font\twelvebsy=cmbsy10 at 12pt
\font\tenbsy=cmbsy10
\font\eightbsy=cmbsy8
\font\sevenbsy=cmbsy7
\font\sixbsy=cmbsy6
\font\fivebsy=cmbsy5

\font\tenmsa=msam10

\font\sevenmsa=msam7
\font\fivemsa=msam5
\newfam\msafam
  \textfont\msafam=\tenmsa
  \scriptfont\msafam=\sevenmsa
  \scriptscriptfont\msafam=\fivemsa

\font\tenmsb=msbm10
\font\eightmsb=msbm8
\font\sevenmsb=msbm7
\font\fivemsb=msbm5
\newfam\msbfam
  \textfont\msbfam=\tenmsb
  \scriptfont\msbfam=\sevenmsb
  \scriptscriptfont\msbfam=\fivemsb
\def\Bbb{\fam\msbfam\tenmsb}

\font\tenCal=eusm10
\font\sevenCal=eusm7
\font\fiveCal=eusm5
\newfam\Calfam
  \textfont\Calfam=\tenCal
  \scriptfont\Calfam=\sevenCal
  \scriptscriptfont\Calfam=\fiveCal
\def\Cal{\fam\Calfam\tenCal}

\font\teneuf=eusm10
\font\teneuf=eufm10
\font\seveneuf=eufm7
\font\fiveeuf=eufm5
\newfam\euffam
  \textfont\euffam=\teneuf
  \scriptfont\euffam=\seveneuf
  \scriptscriptfont\euffam=\fiveeuf

\font\seventeenbfit=cmmib10 at 17.28pt

\font\twelvebfit=cmmib10 at 12pt
\font\tenbfit=cmmib10
\font\eightbfit=cmmib8
\font\sevenbfit=cmmib7
\font\sixbfit=cmmib6
\font\fivebfit=cmmib5
\newfam\bfitfam
  \textfont\bfitfam=\tenbfit
  \scriptfont\bfitfam=\sevenbfit
  \scriptscriptfont\bfitfam=\fivebfit


\catcode`\@=11
\def\eightpoint{%
  \textfont0=\eightrm \scriptfont0=\sixrm \scriptscriptfont0=\fiverm
  \def\rm{\fam\z@\eightrm}%
  \textfont1=\eighti \scriptfont1=\sixi \scriptscriptfont1=\fivei
  \def\oldstyle{\fam\@ne\eighti}%
  \textfont2=\eightsy \scriptfont2=\sixsy \scriptscriptfont2=\fivesy
  \textfont\itfam=\eightit
  \def\it{\fam\itfam\eightit}%
  \textfont\slfam=\eightsl
  \def\sl{\fam\slfam\eightsl}%
  \textfont\bffam=\eightbf \scriptfont\bffam=\sixbf
  \scriptscriptfont\bffam=\fivebf
  \def\bf{\fam\bffam\eightbf}%
  \textfont\ttfam=\eighttt
  \def\tt{\fam\ttfam\eighttt}%
  \textfont\msbfam=\eightmsb
  \def\Bbb{\fam\msbfam\eightmsb}%
  \abovedisplayskip=9pt plus 2pt minus 6pt
  \abovedisplayshortskip=0pt plus 2pt
  \belowdisplayskip=9pt plus 2pt minus 6pt
  \belowdisplayshortskip=5pt plus 2pt minus 3pt
  \smallskipamount=2pt plus 1pt minus 1pt
  \medskipamount=4pt plus 2pt minus 1pt
  \bigskipamount=9pt plus 3pt minus 3pt
  \normalbaselineskip=9pt
  \setbox\strutbox=\hbox{\vrule height7pt depth2pt width0pt}%
  \let\bigf@ntpc=\eightrm \let\smallf@ntpc=\sixrm
  \normalbaselines\rm}
\catcode`\@=12

\def\eightpointbf{%
 \textfont0=\eightbf   \scriptfont0=\sixbf   \scriptscriptfont0=\fivebf
 \textfont1=\eightbfit \scriptfont1=\sixbfit \scriptscriptfont1=\fivebfit
 \textfont2=\eightbsy  \scriptfont2=\sixbsy  \scriptscriptfont2=\fivebsy
 \eightbf
 \baselineskip=10pt}

\def\tenpointbf{%
 \textfont0=\tenbf   \scriptfont0=\sevenbf   \scriptscriptfont0=\fivebf
 \textfont1=\tenbfit \scriptfont1=\sevenbfit \scriptscriptfont1=\fivebfit
 \textfont2=\tenbsy  \scriptfont2=\sevenbsy  \scriptscriptfont2=\fivebsy
 \tenbf}

\def\twelvepointbf{%
 \textfont0=\twelvebf   \scriptfont0=\eightbf   \scriptscriptfont0=\sixbf
 \textfont1=\twelvebfit \scriptfont1=\eightbfit \scriptscriptfont1=\sixbfit
 \textfont2=\twelvebsy  \scriptfont2=\eightbsy  \scriptscriptfont2=\sixbsy
 \twelvebf
 \baselineskip=14.4pt}

\def\seventeenpointbf{%
 \textfont0=\seventeenbf  \scriptfont0=\twelvebf  \scriptscriptfont0=\eightbf
 \textfont1=\seventeenbfit\scriptfont1=\twelvebfit\scriptscriptfont1=\eightbfit
 \textfont2=\seventeenbsy \scriptfont2=\twelvebsy \scriptscriptfont2=\eightbsy
 \seventeenbf
 \baselineskip=20.736pt}


\newdimen\srdim \srdim=\hsize
\newdimen\irdim \irdim=\hsize
\def\NOSECTREF#1{\noindent\hbox to \srdim{\null\dotfill ???(#1)}}
\def\SECTREF#1{\noindent\hbox to \srdim{\csname REF\romannumeral#1\endcsname}}
\def\INDREF#1{\noindent\hbox to \irdim{\csname IND\romannumeral#1\endcsname}}
\newlinechar=`\^^J
\def\openauxfile{
  \immediate\openin1\jobname.aux
  \ifeof1
  \message{^^JCAUTION\string: you MUST run TeX a second time^^J}
  \let\sectref=\NOSECTREF \let\indref=\NOSECTREF
  \else
  \input \jobname.aux
  \message{^^JCAUTION\string: if the file has just been modified you may
    have to run TeX twice^^J}
  \let\sectref=\SECTREF \let\indref=\INDREF
  \fi
  \message{to get correct page numbers displayed in Contents or Index
    Tables^^J}
  \immediate\openout1=\jobname.aux
  \let\END=\end \def\end{\immediate\closeout1\END}}

\newbox\titlebox   \setbox\titlebox\hbox{\hfil}
\newbox\sectionbox \setbox\sectionbox\hbox{\hfil}
\def\folio{\ifnum\pageno=1 \hfil \else \ifodd\pageno
           \hfil {\eightpoint\copy\sectionbox\kern8mm\number\pageno}\else
           {\eightpoint\number\pageno\kern8mm\copy\titlebox}\hfil \fi\fi}
\footline={\hfil}
\headline={\folio}

\def\titlerunning#1{\setbox\titlebox\hbox{\eightpoint #1}}
\def\title#1{\noindent\hfil$\smash{\hbox{\seventeenpointbf #1}}$\hfil
             \titlerunning{#1}\medskip}

\newcount\numbersection \numbersection=-1
\def\sectionrunning#1{\setbox\sectionbox\hbox{\eightpoint #1}
  \immediate\write1{\string\def \string\REF
      \romannumeral\numbersection \string{%
      \noexpand#1 \string\dotfill \space \number\pageno \string}}}
\def\section#1{%
  \par\vskip0.666cm\penalty -100
  \vbox{\baselineskip=14.4pt\noindent{{\twelvepointbf #1}}}
  \vskip2pt
  \penalty 500
  \advance\numbersection by 1
  \sectionrunning{#1}}

\def\subsection#1{%
  \par\vskip0.5cm\penalty -100
  \vbox{\noindent{{\tenpointbf #1}}}
  \vskip1pt
  \penalty 500}

\newcount\numberindex \numberindex=0
\def\index#1#2{%
  \advance\numberindex by 1
  \immediate\write1{\string\def \string\IND #1%
     \romannumeral\numberindex \string{%
     \noexpand#2 \string\dotfill \space \string\S \number\numbersection,
     p.\string\ \space\number\pageno \string}}}

\newdimen\itemindent \itemindent=\parindent

\def\item#1{\par\noindent\hangindent\itemindent%
            \rlap{#1}\kern\itemindent\ignorespaces}
\def\itemitem#1{\par\noindent\hangindent2\itemindent%
            \kern\itemindent\rlap{#1}\kern\itemindent\ignorespaces}
\def\itemitemitem#1{\par\noindent\hangindent3\itemindent%
            \kern2\itemindent\rlap{#1}\kern\itemindent\ignorespaces}

\long\def\claim#1|#2\endclaim{\par\vskip 5pt\noindent
{\tenpointbf #1.}\ {\it #2}\par\vskip 5pt}

\def\proof{\noindent{\it Proof}}

\def\today{\ifcase\month\or
January\or February\or March\or April\or May\or June\or July\or August\or
September\or October\or November\or December\fi \space\number\day,
\number\year}

\catcode`\@=11
\newcount\@tempcnta \newcount\@tempcntb
\def\timeofday{{%
\@tempcnta=\time \divide\@tempcnta by 60 \@tempcntb=\@tempcnta
\multiply\@tempcntb by -60 \advance\@tempcntb by \time
\ifnum\@tempcntb > 9 \number\@tempcnta:\number\@tempcntb
  \else\number\@tempcnta:0\number\@tempcntb\fi}}
\catcode`\@=12

\def\bibitem#1&#2&#3&#4&%
{\hangindent=1.8cm\hangafter=1
\noindent\rlap{\hbox{\eightpointbf #1}}\kern1.8cm{\rm #2}{\it #3}{\rm #4.}}


\def\bC{{\Bbb C}}

\def\bN{{\Bbb N}}
\def\bP{{\Bbb P}}
\def\bQ{{\Bbb Q}}
\def\bR{{\Bbb R}}

\def\bZ{{\Bbb Z}}


\def\cC{{\Cal C}}
\def\cD{{\Cal D}}
\def\cE{{\Cal E}}
\def\cF{{\Cal F}}
\def\cG{{\Cal G}}

\def\cI{{\Cal I}}

\def\cK{{\Cal K}}

\def\cM{{\Cal M}}
\def\cN{{\Cal N}}
\def\cO{{\Cal O}}

\def\cX{{\Cal X}}
\def\la{{\longrightarrow}}


\def\square{{\hfill \hbox{
\vrule height 1.453ex  width 0.093ex  depth 0ex
\vrule height 1.5ex  width 1.3ex  depth -1.407ex\kern-0.1ex
\vrule height 1.453ex  width 0.093ex  depth 0ex\kern-1.35ex
\vrule height 0.093ex  width 1.3ex  depth 0ex}}}
\def\qed{\kern10pt$\square$}
\def\hexnbr#1{\ifnum#1<10 \number#1\else
 \ifnum#1=10 A\else\ifnum#1=11 B\else\ifnum#1=12 C\else
 \ifnum#1=13 D\else\ifnum#1=14 E\else\ifnum#1=15 F\fi\fi\fi\fi\fi\fi\fi}
\def\msatype{\hexnbr\msafam}
\def\msbtype{\hexnbr\msbfam}
\mathchardef\restriction="3\msatype16   
\mathchardef\boxsquare="3\msatype03
\mathchardef\preccurlyeq="3\msatype34
\mathchardef\compact="3\msatype62
\mathchardef\smallsetminus="2\msbtype72   \let\ssm\smallsetminus
\mathchardef\subsetneq="3\msbtype28
\mathchardef\supsetneq="3\msbtype29
\mathchardef\leqslant="3\msatype36   \let\le\leqslant
\mathchardef\geqslant="3\msatype3E   \let\ge\geqslant
\mathchardef\stimes="2\msatype02
\mathchardef\ltimes="2\msbtype6E
\mathchardef\rtimes="2\msbtype6F

\def\smallvee{{\scriptscriptstyle\vee}}

\def\ddbar{\partial\overline\partial}

\let\ol=\overline

\let\wt=\widetilde
\let\wh=\widehat
\let\text=\hbox
\def\buildo#1^#2{\mathop{#1}\limits^{#2}}
\def\buildu#1_#2{\mathop{#1}\limits_{#2}}
\def\ort{\mathop{\hbox{\kern1pt\vrule width4.0pt height0.4pt depth0pt
    \vrule width0.4pt height6.0pt depth0pt\kern3.5pt}}}
\let\lra\longrightarrow
\def\vlra{\mathrel{\smash-}\joinrel\mathrel{\smash-}\joinrel%
    \kern-2pt\longrightarrow}
\def\srelbar{\vrule width0.6ex height0.65ex depth-0.55ex}
\def\merto{\mathrel{\srelbar\kern1.3pt\srelbar\kern1.3pt\srelbar
    \kern1.3pt\srelbar\kern-1ex\raise0.28ex\hbox{${\scriptscriptstyle>}$}}}


\def\Hom{\mathop{\rm Hom}\nolimits}

\def\Pic{\mathop{\rm Pic}\nolimits}
\def\Def{\mathop{\rm Def}\nolimits}
\def\num{\nu}

\def\rank{\mathop{\rm rank}\nolimits}

\def\nef{\mathop{\rm nef}\nolimits}

\def\codim{\mathop{\rm codim}\nolimits}

\def\Sing{\mathop{\rm Sing}\nolimits}

\def\NE{\mathop{\rm NE}\nolimits}
\def\ME{\mathop{\rm ME}\nolimits}
\def\SME{\mathop{\rm SME}\nolimits}
\def\Vol{\mathop{\rm Vol}\nolimits}

\def\NS{{\rm NS}}
\def\FS{{\rm FS}}
\def\nnef{{\rm nnef}}

\long\def\InsertFig#1 #2 #3 #4\EndFig{\par
\hbox{\hskip #1mm$\vbox to#2mm{\vfil\special{"
(/home/demailly/psinputs/grlib.ps) run
#3}}#4$}}
\long\def\LabelTeX#1 #2 #3\ELTX{\rlap{\kern#1mm\raise#2mm\hbox{#3}}}


\itemindent = 7mm

\title{The pseudo-effective cone of a} 
\title{compact K\"ahler manifold and}
\title{varieties of negative Kodaira dimension}
\titlerunning{The pseudo-effective cone of compact K\"ahler manifolds}
\vskip10pt

{\noindent\hangindent0.6cm\hangafter-2
{\bf S\'ebastien Boucksom${}^1$ \hfill Jean-Pierre Demailly${}^2$\kern0.45cm\\
     Mihai Paun${}^3$ \hfill  Thomas Peternell${}^4$\kern0.45cm

}}

{\noindent\hangindent0.6cm\hangafter-4{\it
\llap{${}^1$}Universit\'e de Paris VII\hfill
\llap{${}^2$}Universit\'e de Grenoble I, BP 74\kern0.6cm\break
Institut de Math\'ematiques\hfill Institut Fourier\kern0.6cm\break
175 rue du Chevaleret\hfill UMR 5582 du CNRS\kern0.6cm\break
75013 Paris, France\hfill 38402 Saint-Martin d'H\`eres, France\kern0.6cm
}
\vskip2pt
{\noindent\hangindent0.6cm\hangafter-4{\it
\llap{${}^3$}Universit\'e de Strasbourg\kern0.6cm\hfill
\llap{${}^4$}Universit\"at Bayreuth\kern0.6cm\break
D\'epartement de Math\'ematiques\hfill Mathematisches Institut\kern0.6cm\break
67084 Strasbourg, France\hfill D-95440 Bayreuth, Deutschland\kern0.6cm

}}
\vskip20pt

\noindent{\bf Abstract. \it { We prove that a holomorphic line bundle 
    on a projective manifold is pseudo-effective if and only if its
    degree on any member of a covering family of curves is non-negative. 
    This is a consequence of a duality statement between the cone of
    pseudo-effective divisors and the cone of ``movable curves'', which
    is obtained from a general theory of movable intersections and
    approximate Zariski decomposition for closed positive $(1,1)$-currents.
    As a corollary, a projective manifold has a pseudo-effective canonical
    bundle if and only if it is is not uniruled.  We also prove that a
    4-fold with a canonical bundle which is pseudo-effective and of numerical 
    class zero in restriction to curves of a covering family, has non 
    negative Kodaira dimension. }  }
\vskip20pt

\section{\S0 Introduction}

One of the major open problems in the classification theory of
projective or compact K\"ahler manifolds is the following geometric
description of varieties of negative Kodaira dimension.

\claim 0.1 Conjecture| A projective $($or compact K\"ahler$)$ manifold $X$
has Kodaira dimension $\kappa(X) = - \infty$ if and only if $X$ is
uniruled.
\endclaim

One direction is trivial, namely $X$ uniruled implies $\kappa(X)=-\infty$.
Also, the conjecture is known to be true for projective threefolds
by [Mo88] and for non-algebraic K\"ahler threefolds by [Pe01], with the 
possible exception of simple threefolds (recall that a variety is said to 
be simple if there is no compact positive dimensional subvariety
through a very general point of~$X$). In the case of projective
manifolds, the problem can be split into two more tractable parts~:

\medskip
{\itemindent 5mm
\item {A.} {\it If the canonical bundle $K_X$ is not pseudo-effective, 
i.e.\ not contained in the closure of the cone spanned by classes of 
effective divisors, then $X$ is uniruled.}
\item {B.} {\it If $K_X$ is pseudo-effective, then $\kappa(X) \geq 0$. }
}

\medskip 
In the K\"ahler case, the statements should be essentially the same, except
that effective divisors have to be replaced by closed positive 
$(1,1)$-currents.

In this paper we give a positive answer to (A) for projective
manifolds of any dimension, and a partial answer to (B) for $4$-folds.
Part (A) follows in fact from a much more general fact which describes 
the geometry of the pseudo-effective cone.

\claim 0.2 Theorem| A line bundle $L$ on a projective manifold $X$ is
pseudo-effective if and only if $L \cdot C \geq 0 $ for all
irreducible curves $C$ which move in a family covering $X$.
\endclaim

In other words, the dual cone to the pseudo-effective cone is the
closure of the cone of ``movable'' curves. This should be compared with the
duality between the nef cone and the cone of effective curves.

\claim 0.3 Corollary {\rm (Solution of (A))}|Let $X$ be a projective
manifold. If $K_X$ is not pseudo-effective, then $X$ is covered by
rational curves.
\endclaim

In fact, if $K_X$ is not pseudo-effective, then by (0.2) there exists
a covering family $(C_t)$ of curves with $K_X \cdot C_t < 0$, so that
(0.3) follows by a well-known characteristic $p$ argument of Miyaoka 
and Mori [MM86] (the so called bend-and-break lemma essentially amounts
to deform the $C_t$ so that they break into pieces, one of which is
a rational curve).  

In the K\"ahler case both a suitable analogue to (0.2) and the theorem of
Miyaoka-Mori are unknown. It should also be mentioned that the duality
statement following (0.2) is actually (0.2) for $\bR$-divisors. The
proof is based on a use of ``approximate Zariski decompositions'' and
an estimate for an intersection number related to this decomposition.
A major tool is the volume of an $\bR$-divisor which distinguishes big
divisors (positive volume) from divisors on the boundary of the
pseudo-effective cone (volume $0$).

Concerning (B) we prove the following weaker statement.

\claim 0.4 Theorem| Let $X$ be a
smooth projective $4$-fold. Assume that $K_X$ is pseudo-effective and
there is a covering family $(C_t)$ of curves such that $K_X \cdot C_t = 0$. 
Then $\kappa(X) \geq 0$.
\endclaim

One important ingredient of the proof of (0.4) is the quotient defined
by the family $(C_t)$. In order to obtain the full answer to Problem
(B) in dimension $4$, we would still need to prove that $K_X$ is
effective if $K_X$ is positive on all covering families of curves. In
fact, in that case, $K_X$ should be big, i.e.\ of maximal Kodaira
dimension.

\section{\S1 Positive cones in the spaces of divisors and of curves}

In this section we introduce the relevant cones, both in the
projective and K\"ahler contexts -- in the latter case, divisors and curves
should simply be replaced by positive currents of bidimension $(n-1,n-1)$
and $(1,1)$, respectively. We implicitly use that all (De Rham,
resp.\ Dolbeault) cohomology groups under consideration can be
computed in terms of smooth forms or currents, since in both cases we
get resolutions of the same sheaf of locally constant functions
(resp.\ of holomorphic sections).

\claim 1.1 Definition| Let $X$ be a compact K\"ahler manifold.
\smallskip
\item{\rm(i)} The K\"ahler cone is the set 
$\cK\subset H^{1,1}_\bR(X)$ of classes $\{\omega\}$
of K\"ahler forms $($this is an open convex cone$)$.
\smallskip
\item{\rm(ii)} The pseudo-effective cone is the set
$\cE\subset H^{1,1}_\bR(X)$ of classes $\{T\}$ of closed positive 
currents of type $(1,1)$ $($this is a closed convex cone$)$.
Clearly $\cE\supset\overline{\cK}$.
\smallskip
\item {\rm (iii)} The Neron-Severi
space is defined by
$$
\NS_\bR(X):= \big(H^{1,1}_\bR(X)\cap H^2(X,\bZ)/{\rm tors}\big)\otimes_\bZ\bR.
$$ 
\item {\rm (iv)} We set
$$\cK_\NS = \cK\cap NS_\bR(X),\qquad \cE_\NS=\cE\cap NS_{\bR}(X).$$
\endclaim

Algebraic geometers tend to restrict themselves to the algebraic cones
generated by ample divisors and effective divisors, respectively.
Using $L^2$ estimates for $\overline\partial$, one can show the following
expected relations between the algebraic and transcendental cones (see
[De90], [De92]).

\claim 1.2 Proposition|In a projective manifold~$X$,
$\cE_\NS$ is the closure of the convex cone generated by effective
divisors, and $\overline{\cK_\NS}$ is the closure of the cone 
generated by nef $\bR$-divisors.
\endclaim

By extension, we will say that $\overline{\cK}$ is the cone of {\it nef}
$(1,1)$-cohomology classes (even though they are not necessarily integral).
We now turn ourselves to cones in cohomo\-logy of bidegree $(n-1,n-1)$.
\bigskip

\claim 1.3 Definition|
Let $X$ be a compact K\"ahler manifold.
\smallskip
\item {\rm (i)} We define $\cN$ to be the $($closed$)$ convex cone in 
$H^{n-1,n-1}_\bR(X)$ gene\-rated by classes of positive currents $T$ of type
\hbox{$(n-1,n-1)$} $($i.e., of bidimension $(1,1))$.

\item{\rm (ii)} We define the cone $\cM\subset H^{n-1,n-1}_\bR(X)$ of
{\bf movable classes} to be the closure of the convex cone generated by
classes of currents of the form
$$\mu_\star(\wt\omega_1\wedge\ldots\wedge\wt\omega_{n-1})$$
where $\mu:\wt X\to X$ is an arbitrary modification $($one could just
restrict oneself to compositions of blow-ups with smooth centers$)$, and 
the $\wt\omega_j$ are K\"ahler forms on~$\wt X$. Clearly $\cM\subset\cN$.

\item {\rm (iii)} Correspondingly, we introduce the intersections
$$ 
\cN_\NS=\cN\cap N_1(X),\qquad \cM_\NS=\cM\cap N_1(X),\qquad 
$$
in the space of integral bidimension $(1,1)$-classes
$$
N_1(X):=(H^{n-1,n-1}_\bR(X)\cap H^{2n-2}(X,\bZ)/{\rm tors})\otimes_\bZ\bR.
$$
\item{\rm (iv)} If $X$ is projective, we define $\NE(X)$ to be the
convex cone generated by all effective curves. Clearly  
$\overline{\NE(X)}\subset\cN_\NS$.

\item{\rm (v)} If $X$ is projective, we say that $C$ is a  {\bf 
strongly movable curve} if 
$$
C=\mu_\star(\wt A_1\cap\ldots\cap \wt A_{n-1})
$$ 
for suitable very ample divisors $\wt A_j$ on $\wt X$, where 
$\mu:\wt X\to X$ is a modification. We let ${\rm SME}(X)$ 
to be the convex cone generated by all strongly movable $($effective$)$ 
curves. Clearly  $\overline{{\rm SME}(X)}\subset\cM_\NS$.
\smallskip
\item{\rm (vi)} We say that $C$ is a { \bf movable curve} if $C=C_{t_0}$
is a member of an analytic family $(C_t)_{t\in S}$ such that 
$\bigcup_{t\in S}C_t=X$
and, as such, is a reduced irreducible $1$-cycle. We let ${\rm ME}(X)$ 
to be the convex cone generated by all movable $($effective$)$ curves.
\vskip0pt
\endclaim

The upshot of this definition lies in the following easy observation.

\claim 1.4 Proposition|Let $X$ be a compact K\"ahler manifold. 
Consider the Poincar\'e duality pairing
$$
H^{1,1}_\bR(X)\times H^{n-1,n-1}_\bR(X)\lra \bR,\qquad
(\alpha,\beta)\longmapsto\int_X \alpha\wedge\beta.
$$
Then the duality pairing takes nonnegative values
\smallskip
\item{\rm (i)} for all pairs $(\alpha,\beta)\in\ol\cK\times\cN$;
\item{\rm (ii)} for all pairs $(\alpha,\beta)\in\cE\times\cM$.
\item{\rm (iii)} for all pairs $(\alpha,\beta)$ where $\alpha\in\cE$ and 
$\beta=[C_t]\in{\rm ME}(X)$ is the class of a movable curve.
\vskip0pt
\endclaim

\proof. (i) is obvious. In order to prove (ii), we may assume that
$\beta=\mu_\star(\wt\omega_1\wedge\ldots\wedge\wt\omega_{n-1})$ for
some modification $\mu:\wt X\to X$, where $\alpha=\{T\}$ is the class
of a positive $(1,1)$-current on $X$ and $\wt\omega_j$ are K\"ahler
forms on~$\wt X$. Then
$$
\int_X\alpha\wedge\beta=\int_X T\wedge
\mu_\star(\wt\omega_1\wedge\ldots\wedge\wt\omega_{n-1})=
\int_X \mu^* T\wedge
\wt\omega_1\wedge\ldots\wedge\wt\omega_{n-1}\ge 0.
$$
Here, we have used the fact that a closed positive $(1,1)$-current $T$
always has a pull-back $\mu^\star T$, which follows from the fact that
if $T=i\ddbar\varphi$ locally for some plurisubharmonic function in~$X$, 
we can set 
$\mu^\star T=i\ddbar(\varphi\circ\mu)$. For (iii), we suppose
$\alpha=\{T\}$ and $\beta=\{[C_t]\}$. Then we take an open covering
$(U_j)$ on $X$ such that $T=i\ddbar\varphi_j$ with suitable plurisubharmonic 
functions $\varphi_j$ on $U_j$. If we select a smooth partition of unity 
$\sum\theta_j=1$ subordinate to $(U_j)$, we then get
$$
\int_X\alpha\wedge\beta=\int_{C_t} T_{|C_t}=\sum_j\int_{C_t\cap U_j}
\theta_j i\ddbar\varphi_{j|C_t}\ge 0.
$$
For this to make sense, it should be noticed that $T_{|C_t}$ is a well 
defined closed positive $(1,1)$-current (i.e.\ measure) on $C_t$ for almost 
every $t\in S$, in the sense of Lebesgue measure. This is true only because
$(C_t)$ covers $X$, thus $\varphi_{j|C_t}$ is not identically
$-\infty$ for almost every $t\in S$. The equality in the last formula
is then shown by a regularization argument for $T$, writing $T=\lim T_k$
with $T_k=\alpha+i\ddbar\psi_k$ and a decreasing sequence of smooth almost 
plurisubharmonic potentials $\psi_k\downarrow\psi$ such that the Levi forms 
have a uniform lower bound $i\ddbar\psi_k\ge -C\omega$ (such a
sequence exists by [De92]). Then, writing $\alpha=i\ddbar v_j$ for some
smooth potential $v_j$ on $U_j$, we have $T=i\ddbar\varphi_j$ on $U_j$
with $\varphi_j=v_j+\psi$, and this is the decreasing limit of the smooth
approximations $\varphi_{j,k}=v_j+\psi_k$ on $U_j$. Hence 
$T_{k|C_t}\to T_{|C_t}$ for the weak topology of measures on $C_t$.\qed
\bigskip

If $\cC$ is a convex cone in a finite dimensional vector space $E$,
we denote by $\cC^\smallvee$ the dual cone, i.e.\ the set of linear
forms $u\in E^\star$ which take nonnegative values on all elements of $\cC$. 
By the Hahn-Banach theorem, we always have $\cC^{\smallvee\smallvee}=
\overline{\cC}$. 

Proposition 1.4 leads to the natural question whether the 
cones $(\cK,\cN)$ and $(\cE,\cM)$ are dual under Poincar\'e duality. 
This question is addressed in the next section. Before doing so,
we observe that the algebraic and transcendental cones
of $(n-1,n-1)$ cohomology classes are related by the following 
equalities (similar to what we already noticed for $(1,1)$-classes, 
see Prop.~1.2).

\claim 1.5 Theorem|Let $X$ be a projective manifold. Then
\smallskip
\item{\rm(i)} $\overline{\NE(X)} =\cN_\NS$.
\smallskip
\item{\rm(ii)} $\overline{{\rm SME}(X)}=\overline{{\rm ME}(X)} =\cM_\NS$.
\endclaim

\proof. (i) It is a standard result of algebraic geometry (see e.g.\
[Ha70]), that the cone of effective cone $\NE(X)$ is dual to the cone
$\overline{\cK_\NS}$ of nef divisors, hence 
$$\cN_\NS\supset\overline{\NE(X)}=\cK^\vee.$$ 
On the other hand, (1.4)~(i) implies that $\cN_\NS\subset\cK^\vee$, so we must
have equality and (i) follows.

Similarly, (ii) requires a duality statement which will be established only
in the next sections, so we postpone the proof.\qed

\section{\S2 Main results and conjectures}

First, the already mentioned duality between nef divisors and effective 
curves extends to the K\"ahler case and to transcendental classes. More 
precisely, [DPa03] gives

\claim 2.1 Theorem {\rm(Demailly-Paun, 2001)}|If $X$ is K\"ahler, then 
the cones $\overline{\cK}\subset H^{1,1}_\bR(X)$ and 
$\cN\subset H^{n-1,n-1}_\bR(X)$ are dual by Poincar\'e duality, and
$\cN$ is the closed convex cone generated by classes $[Y]\wedge\omega^{p-1}$
where $Y\subset X$ ranges over $p$-dimensional analytic subsets,
$p=1,2,\ldots, n$, and $\omega$ ranges over K\"ahler forms.
\endclaim

\proof. Indeed, Prop.~1.4 shows that the dual cone $\cK^\vee$ contains
$\cN$ which itself contains the cone $\cN'$ of all classes of the form
$\{[Y]\wedge\omega^{p-1}\}$. The main result of [DPa03] conversely shows
that the dual of $(\cN')^\vee$ is equal to $\overline\cK$, so we must
have
$$\cK^\vee = \overline{\cN'} = \cN.\eqno\square$$

The main new result of this paper is the following characterization of
pseudo-effective classes (in which the ``only if'' part already follows 
from 1.4~(iii)).

\claim 2.2 Theorem|If $X$ is projective, then a class $\alpha\in\NS_\bR(X)$
is pseudo-effective if $($and only if$\,)$ it is in the dual cone of the cone
$\SME(X)$ of strongly movable curves.
\endclaim

In other words, a line bundle $L$ is pseudo-effective if (and
only if) $L\cdot C\ge 0$ for all {\it movable curves}, i.e., $L\cdot C\ge 0$
for every very generic curve $C$ (not contained in a countable union of 
algebraic subvarieties). In fact, by definition of $\SME(X)$, it is enough 
to consider only those curves $C$ which are images of generic complete 
intersection of very ample divisors on some variety $\wt X$, under
a modification $\mu:\wt X\to X$.

By a standard blowing-up argument, it also follows that a line bundle
$L$ on a normal Moishezon variety is pseudo-effective if and only if
$L \cdot C \geq 0$ for every movable curve~$C$. 
\smallskip

The K\"ahler analogue should be~:

\claim 2.3 Conjecture| For an arbitrary compact 
K\"ahler manifold $X$, the cones $\cE$ and $\cM$ are dual.
\endclaim

The relation between the various cones of movable curves and currents in 
(1.5) is now a rather direct consequence of Theorem 2.2. In fact, using 
ideas hinted in [DPS96], we can say a little bit more. Given 
an irreducible curve $C \subset X$, we consider its normal ``bundle''  
$N_C = \Hom(\cI/\cI^2,\cO_C)$, where $\cI$ is the ideal sheaf of $C$. 
If $C$ is a general member of a covering family $(C_t)$, then $N_C$ 
is nef. Now [DPS96] says that the dual cone of the pseudo-effective cone
of $X$ contains the closed cone spanned by curves with nef normal 
bundle, which in turn contains the cone of movable curves. In this
way we get~:

\claim 2.4 Theorem|Let $X$ be a projective manifold. Then the following 
cones coincide.
\smallskip
\item {\rm (i)} the cone $\cM_\NS=\cM\cap N_1(X)\,;$
\smallskip
\item {\rm (ii)} the closed cone $\overline{\SME(X)}$ of strongly movable 
curves$\,;$
\smallskip
\item {\rm (iii)} the closed cone $\overline{\ME(X)}$ of movable curves$\,;$
\smallskip
\item {\rm (iv)} the closed cone $\overline{\ME_{\nef}(X)}$ of curves with nef 
normal bundle.\vskip0pt
\endclaim

\proof. We have already seen that 
$$\SME(X)\subset\ME(X)\subset\ME_{\nef}(X)\subset(\cE_\NS)^\smallvee$$
and 
$$
\SME(X)\subset\ME(X)\subset\cM_\NS\subset(\cE_\NS)^\smallvee
$$ 
by 1.4~(iii). Now Theorem 2.2 implies $(\cM_\NS)^\smallvee =
\overline{\SME(X)}$, and 2.4 follows.\qed

\claim 2.5 Corollary| Let $X$ be a projective manifold and $L$ a line bundle 
on $X$.
\smallskip
\item {\rm (i)} $L$ is pseudo-effective if and only if $L \cdot C \geq 0$ 
for all curves $C$ with nef
normal sheaf $N_C$.
\smallskip
\item {\rm (ii)} If $L$ is big, then $L \cdot C > 0$ for all curves $C$ 
with nef normal sheaf $N_C$.\vskip0pt
\endclaim

2.5~(i) strenghtens results from [PSS99]. It is however not yet clear whether 
$\cM_\NS=\cM\cap N_1(X)$ is equal to the closed cone of curves with 
{\it ample} normal bundle (although we certainly expect this to be true).

The most important special case of Theorem 2.2 is

\claim 2.6 Theorem| If $X$ is a projective manifold and is not uniruled,
then $K_X$ is pseudo-effective, i.e.\ $K_X\in\cE_\NS$.
\endclaim

\proof. This is merely a restatement of Corollary~0.3, which was proved in
the introduction (as a consequence of the results of [MM86]). 
\smallskip 

Theorem 2.6 can be generalized as follows.

\claim 2.7 Theorem| Let $X$ be a projective manifold $($or a normal 
projective variety$\,)$. Let $\cF\subset T_X$ be a coherent subsheaf.
If $\det \cF^*$ is not pseudo-effective, then $X$ is uniruled. 
In other words, if $X$ is not uniruled and $\Omega^1_X \to\cG$ is 
generically surjective, then $\det \cG$ is pseudo-effective. 
\endclaim 

\proof. In fact, since $\det \cF^*$ is not pseudo-effective, there exists by
(2.2) a covering family $(C_t)$ such that $c_1(\cF) \cdot C_t > 0$.
Hence $X$ is uniruled by [Mi87], [SB92]. \qed

\claim 2.8. Problems| 
\item {\rm(1)} Does 2.7 generalize to subsheaves
$\cF \subset T_X^{\otimes m}\,?$
\item {\rm(2)} Suppose in 2.7 that only $\kappa(\det \cF^*) = - \infty$. 
Is $X$ still uniruled$\,?$ 
What can be said if $c_1(\cF^*)$ is on the boundary of the 
pseudo-effective cone$\,?$\vskip0pt
\endclaim 

Turning to varieties with pseudo-effective canonical bundles, we have the

\claim 2.9 Conjecture {\rm (special case of the ``abundance conjecture'')}|
If $K_X$ is pseudo-effective, then $\kappa(X)\ge 0$.
\endclaim

In the last section we will prove this in dimension 4 under the additional assumption
that there is a covering family of curves $(C_t)$ such that $K_X \cdot C_t = 0$. 

\section{\S3 Zariski decomposition and movable intersections}

Let $X$ be compact K\"ahler and let $\alpha\in\cE^\circ$ be in the
{\it interior} of the pseudo--effective cone. In analogy with the
algebraic context such a class $\alpha$ is called ``big'', and it can
then be represented by a {\it K\"ahler current} $T$, i.e.\ a closed
positive $(1,1)$-current $T$ such that $T\ge \delta\omega$ for some
smooth hermitian metric $\omega$ and a constant $\delta\ll 1$.

\claim 3.1 Theorem {\rm (Demailly [De92], [Bou02b, 3.1.24]}|If $T$
is a K\"ahler current, then one can write $T=\lim T_m$ for a sequence of 
K\"ahler currents $T_m$ which have
logarithmic poles with coefficients in ${1\over m}\bZ$, i.e.\ there
are modifications $\mu_m:X_m\to X$ such that
$$\mu_m^\star T_m=[E_m]+\beta_m$$
where $E_m$ is an effective $\bQ$-divisor on $X_m$ with coefficients
in ${1\over m}\bZ$ $($the ``fixed part''$)$ and $\beta_m$ is a closed 
semi-positive form $($the ``movable part''$)$.
\endclaim

\proof. Since this result has already been studied extensively, we just recall
the main idea. Locally we can write $T=i\ddbar\varphi$ for some strictly 
plurisubharmonic potential $\varphi$. By a Bergman kernel trick and the
Ohsawa-Takegoshi $L^2$ extension theorem, we get local approximations
$$
\varphi=\lim\varphi_m,\qquad \varphi_m(z)={1\over 2m}\log\sum_\ell
|g_{\ell,m}(z)|^2
$$
where $(g_{\ell,m})$ is a Hilbert basis of the space of holomorphic
functions which are $L^2$ with respect to the weight $e^{-2m\varphi}$.
This Hilbert basis is also a family of local generators of the globally defined
multiplier ideal sheaf $\cI(mT)=\cI(m\varphi)$. Then $\mu_m:X_m\to X$
is obtained by blowing-up this ideal sheaf, so that 
$$\mu_m^\star\cI(mT)=\cO(-mE_m).$$
We should notice that by approximating $T-{1\over m}\omega$ instead of $T$,
we can replace $\beta_m$ by $\beta_m+{1\over m}\mu^\star\omega$ which is
a big class on $X_m\,$; by playing with the multiplicities of the components
of the exceptional divisor, we could even achieve that $\beta_m$ is a 
K\"ahler class on $X_m$, but this will not be needed here.\qed
\medskip

The more familiar algebraic analogue would be to take \hbox{$\alpha=c_1(L)$}
with a big line bundle $L$ and to blow-up the base locus of $|mL|$,
$m\gg 1$, to get a $\bQ$-divisor decomposition
$$
\mu_m^\star L\sim E_m+D_m,\qquad E_m~~\hbox{effective},~~D_m~~\hbox{free}.
$$
Such a blow-up is usually referred to as a ``log resolution'' of the
linear system $|mL|$, and we say that $E_m+D_m$ is an approximate
Zariski decomposition of $L$. We will also use this terminology for K\"ahler
currents with logarithmic poles.

\claim 3.2 Definition|We define the {\bf volume}, or {\bf movable
self-intersection} of a big class $\alpha\in\cE^\circ$ to be
$$
\Vol(\alpha)=\sup_{T\in \alpha}\int_{\wt X}\beta^n>0
$$
where the supremum is taken over all K\"ahler currents $T\in \alpha$
with logarithmic poles, and $\mu^\star T=[E]+\beta$ with
respect to some modification $\mu:\wt X\to X$.
\endclaim

By Fujita [Fuj94] and Demailly-Ein-Lazarsfeld [DEL00], if $L$ 
is a big line bundle, we have
$$
\Vol(c_1(L))=\lim_{m\to+\infty}D_m^n=\lim_{m\to+\infty}{n!\over m^n}h^0(X,mL),
$$
and in these terms, we get the following statement.

\claim 3.3 Proposition| Let $L$ be a big line bundle on the
projective manifold $X$.  Let $\epsilon > 0$. Then there exists a
modification $\mu: X_{\epsilon} \to X$ and a decomposition $\mu^*(L) =
E + \beta $ with $E$ an effective $\bQ$-divisor and $\beta$ a big and
nef $\bQ$-divisor such that
$$\Vol(L) -\varepsilon\le \Vol(\beta) \le \Vol(L).$$
\endclaim

It is very useful to observe that the supremum in Definition 3.2 can actually
be computed by a collection of currents whose singularities satisfy a 
filtering property.
Namely, if $T_1=\alpha+i\ddbar\varphi_1$ and $T_2=\alpha+i\ddbar\varphi_2$ 
are two K\"ahler currents with logarithmic poles in the class of $\alpha$, 
then 
$$T=\alpha+i\ddbar\varphi,\qquad
\varphi=\max(\varphi_1,\varphi_2)\leqno(3.4)$$ 
is again a K\"ahler current with weaker 
singularities than $T_1$ and $T_2$. One could define as well 
$$T=\alpha+i\ddbar\varphi,\qquad
\varphi={1\over 2m}\log(e^{2m\varphi_1}+e^{2m\varphi_2}),\leqno(3.4')$$
where $m={\rm lcm}(m_1,m_2)$ is the lowest common multiple of the denominators
occuring in $T_1$, $T_2$. Now, take a simultaneous log-resolution
$\mu_m:X_m\to X$ for which the singularities of $T_1$ and $T_2$ are resolved 
as $\bQ$-divisors $E_1$
and $E_2$. Then clearly the associated divisor in the decomposition 
$\mu_m^\star T=[E]+\beta$ is given by $E=\min(E_1,E_2)$. By doing so, 
the volume $\int_{X_m}\beta^n$ gets increased, as we shall see in the
proof of Theorem 3.5 below.

\claim 3.5 Theorem {\rm (Boucksom [Bou02b])}|Let $X$ be a compact K\"ahler
manifold. We denote here by $H^{k,k}_{\ge 0}(X)$ the cone of cohomology classes
of type $(k,k)$ which have non-negative intersection with all closed
semi-positive smooth forms of bidegree $(n-k,n-k)$.
\smallskip
\item{\rm(i)} For each integer $k=1,2,\ldots,n$, there exists a canonical
``movable intersection product''
$$
\cE\times\cdots\times\cE\to H^{k,k}_{\ge 0}(X), \quad
(\alpha_1,\ldots,\alpha_k)\mapsto \langle\alpha_1\cdot\alpha_2\cdots
\alpha_{k-1}\cdot \alpha_k\rangle
$$
such that $\Vol(\alpha)=\langle\alpha^n\rangle$ whenever $\alpha$ is
a big class.
\smallskip
\item{\rm(ii)} The product is increasing, homogeneous of 
degree $1$ and superadditive in each argument, i.e.\
$$
\langle\alpha_1\cdots(\alpha'_j+\alpha''_j)\cdots \alpha_k\rangle\ge
\langle\alpha_1\cdots\alpha'_j\cdots \alpha_k\rangle+
\langle\alpha_1\cdots\alpha''_j\cdots \alpha_k\rangle.
$$
It coincides with the ordinary intersection
product when the $\alpha_j\in\overline{\cK}$ are nef classes.
\smallskip
\item{\rm(iii)} 
The movable intersection product satisfies the Teissier-Hovanskii inequalities
$$
\langle\alpha_1\cdot\alpha_2\cdots \alpha_n\rangle\ge
(\langle\alpha_1^n\rangle)^{1/n}\ldots(\langle\alpha_n^n\rangle)^{1/n}\qquad
\hbox{$($with $\langle\alpha_j^n\rangle=\Vol(\alpha_j)\,)$}.
$$
\smallskip
\item{\rm(iv)} For $k=1$, the above ``product'' reduces to a
$($non linear$)$ projection operator 
$$
\cE\to\cE_1,\qquad \alpha\to\langle\alpha\rangle
$$
onto a certain convex subcone $\cE_1$ of $\cE$ such that
$\overline{\cK}\subset\cE_1\subset\cE$. Moreover, there is
a ``divisorial Zariski decomposition''
$$
\alpha=\{N(\alpha)\}+\langle\alpha\rangle
$$
where $N(\alpha)$ is a uniquely defined effective divisor which is
called the ``negative divisorial part'' of $\alpha$. The map
$\alpha\mapsto N(\alpha)$ is homogeneous and subadditive, and 
$N(\alpha)=0$ if and only if $\alpha\in\cE_1$.
\smallskip
\item{\rm(v)} The components of $N(\alpha)$ always consist of
divisors whose cohomology classes are linearly independent,
especially $N(\alpha)$ has at most $\rho=\rank_\bZ\NS(X)$ components.
\vskip0pt
\endclaim

\proof. We essentially repeat the arguments developped in [Bou02b], with
some simplifications arising from the fact that $X$ is supposed to
be K\"ahler from the start.
\smallskip

\noindent
(i) First assume that all classes $\alpha_j$ are big, i.e.\ 
$\alpha_j\in\cE^\circ$. Fix a smooth closed $(n-k,n-k)$ {\it semi-positive}
form $u$ on $X$. We select K\"ahler currents $T_j\in\alpha_j$ with
logarithmic poles, and a simultaneous log-resolution
$\mu:\wt X\to X$ such that 
$$
\mu^\star T_j=[E_j]+\beta_j.
$$
We consider the direct image current
$\mu_\star(\beta_1\wedge\ldots\wedge\beta_k)$ (which is a closed
positive current of bidegree $(k,k)$ on $X$) and the corresponding
integrals
$$
\int_{\wt X}\beta_1\wedge\ldots\wedge\beta_k\wedge\mu^\star u\ge 0.
$$
If we change the representative $T_j$ with another current $T_j'$,
we may always take a simultaneous log-resolution such that
$\mu^\star T_j'=[E_j']+\beta_j'$, and by using $(3.4')$ we can always
assume that $E'_j\le E_j$. Then $D_j=E_j-E'_j$ is an effective
divisor and we find $[E_j]+\beta_j\equiv [E'_j]+\beta'_j$, hence
$\beta'_j\equiv\beta_j+[D_j]$. A substitution in the integral implies
$$
\eqalign{
\int_{\wt X}\beta_1'\wedge\beta_2
&{}\wedge\ldots\wedge\beta_k\wedge\mu^\star u\cr
&=\int_{\wt X}\beta_1\wedge\beta_2\wedge\ldots\wedge\beta_k\wedge\mu^\star u+
\int_{\wt X}[D_1]\wedge\beta_2\wedge\ldots\wedge\beta_k\wedge\mu^\star u\cr
&\ge \int_{\wt X}\beta_1\wedge\beta_2\wedge\ldots\wedge\beta_k\wedge
\mu^\star u.\cr}
$$
Similarly, we can replace successively all forms $\beta_j$ by the $\beta'_j$,
and by doing so, we find
$$
\int_{\wt X}\beta'_1\wedge\beta'_2\wedge\ldots\wedge\beta'_k\wedge\mu^\star u
\ge\int_{\wt X}\beta_1\wedge\beta_2\wedge\ldots\wedge\beta_k\wedge\mu^\star u.
$$
We claim that the closed positive currents $\mu_\star(\beta_1\wedge\ldots
\wedge\beta_k)$ are uniformly bounded in mass. In fact, if $\omega$ is a 
K\"ahler metric in $X$, there exists a constant $C_j\ge 0$
such that $C_j\{\omega\}-\alpha_j$ is a K\"ahler class. Hence
$C_j\omega-T_j\equiv\gamma_j$ for some K\"ahler form $\gamma_j$ on $X$.
By pulling back with $\mu$, we find
$C_j\mu^\star \omega-([E_j]+\beta_j)\equiv\mu^\star \gamma_j$, hence 
$$\beta_j\equiv C_j\mu^\star \omega-([E_j]+\mu^\star \gamma_j).$$
By performing again a substitution in the integrals, we find
$$
\int_{\wt X}\beta_1\wedge\ldots\wedge\beta_k\wedge\mu^\star u
\le C_1\ldots C_k\int_{\wt X}\mu^\star\omega^k\wedge\mu^\star u
= C_1\ldots C_k\int_X\omega^k\wedge u
$$
and this is true especially for $u=\omega^{n-k}$. We can now arrange
that for each of the integrals associated with a countable dense family 
of forms $u$,  the supremum is achieved by a sequence of currents 
$(\mu_m)_\star(\beta_{1,m}\wedge\ldots\wedge\beta_{k,m})$ obtained
as direct images by a suitable sequence of modifications $\mu_m:
\smash{\wt X_m}\to X$.
By extracting a subsequence, we can achieve that this sequence is weakly 
convergent and we set
$$
\langle\alpha_1\cdot\alpha_2\cdots \alpha_k\rangle =
\mathop{\lim\uparrow}\limits_{m\to+\infty}
\{(\mu_m)_\star(\beta_{1,m}\wedge\beta_{2,m}\wedge\ldots\wedge\beta_{k,m})\}
$$
(the monotonicity is not in terms of the currents themselves, but in terms
of the integrals obtained when we evaluate against a smooth closed
semi-positive form~$u$). By evaluating against a basis of positive classes
$\{u\}\in H^{n-k,n-k}(X)$, we infer by Poincar\'e duality that the class of
$\langle\alpha_1\cdot\alpha_2\cdots \alpha_k\rangle$ is uniquely defined
(although, in general, the representing current is not unique). 
\smallskip

\noindent
(ii) It is indeed clear from the definition that the movable intersection 
product is homogeneous, increasing and superadditive in each argument, at least
when the $\alpha_j$'s are in $\cE^\circ$. However, we can extend the product 
to the closed cone $\cE$ by monotonicity, by setting
$$
\langle\alpha_1\cdot\alpha_2\cdots \alpha_k\rangle =
\mathop{\lim\downarrow}\limits_{\delta\downarrow 0}
\langle(\alpha_1+\delta\omega)\cdot(\alpha_2+\delta\omega)\cdots 
(\alpha_k+\delta\omega)\rangle
$$
for arbitrary classes $\alpha_j\in\cE$ (again, monotonicity occurs
only where we evaluate against closed semi-positive forms~$u$).
By weak compactness, the movable intersection product can always be 
represented by a closed positive current of bidegree $(k,k)$.
\smallskip

\noindent
(iii) The Teissier-Hovanskii inequalities are a direct consequence of the
fact that they hold true for nef classes, so we just have to apply them
to the classes $\beta_{j,m}$ on $\smash{\wt X_m}$ and pass to the limit.
\smallskip

\noindent
(iv) When $k=1$ and $\alpha\in\cE^0$, we have 
$$
\alpha=\lim_{m\to+\infty}\{(\mu_m)_\star T_m\}=
\lim_{m\to+\infty}(\mu_m)_\star[E_m]+\{(\mu_m)_\star\beta_m\}
$$
and $\langle\alpha\rangle=\lim_{m\to+\infty}\{(\mu_m)_\star\beta_m\}$ by 
definition. However, the images $F_m=(\mu_m)_\star E_m$ are effective
$\bQ$-divisors in $X$, and the filtering property implies that $F_m$ is a
decreasing sequence. It must therefore converge to a (uniquely defined)
limit $F=\lim F_m := N(\alpha)$ which is an effective $\bR$-divisor, 
and we get the asserted decomposition in the limit. 

Since $N(\alpha)=\alpha-\langle\alpha\rangle$ we easily see that
$N(\alpha)$ is subadditive and that $N(\alpha)=0$ if $\alpha$ is the
class of a smooth semi-positive form. When $\alpha$ is no longer a big 
class, we define
$$
\langle\alpha\rangle=\lim_{\delta\downarrow 0}\downarrow\langle\alpha+
\delta\omega\rangle,\qquad
N(\alpha)=\lim_{\delta\downarrow 0}\uparrow N(\alpha+\delta\omega)
$$
(the subadditivity of $N$ implies 
$N(\alpha+(\delta+\varepsilon)\omega)\le N(\alpha+\delta\omega)$).
The divisorial Zariski decomposition follows except maybe for the fact 
that $N(\alpha)$ might be a convergent countable sum of divisors. 
However, this will be ruled out when (v) is proved. As 
$N({\scriptscriptstyle\bullet})$ is subadditive
and homogeneous, the set $\cE_1=\{\alpha\in\cE\;;\;N(\alpha)=0\}$ is a 
closed convex cone, and we find that $\alpha\mapsto\langle\alpha\rangle$
is a projection of $\cE$ onto $\cE_1$ (according to [Bou02b], $\cE_1$
consists of those pseudo-effective classes which are ``nef in 
codimension $1$'').
\smallskip

\noindent
(v) Let $\alpha\in\cE^\circ$, and assume that $N(\alpha)$ contains
linearly dependent components $F_j$. Then already all currents
$T\in \alpha$ should be such that $\mu^\star T=[E]+\beta$ where 
$F=\mu_\star E$ contains those linearly dependent components. 
Write $F=\sum\lambda_jF_j$, $\lambda_j>0$ and assume that 
$$\sum_{j\in J}c_jF_j\equiv 0$$
for a certain non trivial linear combination. Then some of the
coefficients $c_j$ must be negative (and some other positive).
Then $E$ is numerically equivalent to
$$
E'\equiv E+t\mu^\star\Big(\sum\lambda_jF_j\Big),
$$
and by choosing $t>0$ appropriate, we obtain an effective divisor $E'$ 
which has a zero coefficient on one of the components $\mu^\star F_{j_0}$.
By replacing $E$ with $\min(E,E')$ via $(3.4')$, we eliminate the component
$\mu^\star F_{j_0}$. This is a contradiction since $N(\alpha)$
was supposed to contain $F_{j_0}$.\qed
\medskip

\claim 3.6 Definition|For a class $\alpha\in H^{1,1}_\bR(X)$, 
we define the {\bf numerical dimension} $\num(\alpha)$ to be 
$\num(\alpha)=-\infty$ if $\alpha$ is not pseudo-effective, and
$$
\num(\alpha)=\max\{p\in\bN\;;\;\langle\alpha^p\rangle\ne 0\},\qquad
\num(\alpha)\in\{0,1,\ldots,n\}
$$
if $\alpha$ is pseudo-effective.
\endclaim

By the results of [DP03], a class is big ($\alpha\in\cE^\circ$) if and only
if $\num(\alpha)=n$. Classes of numerical dimension $0$ can be
described much more precisely, again following Boucksom [Bou02b].

\claim 3.7 Theorem|Let $X$ be a compact
K\"ahler manifold. Then the subset $\cD_0$ of irreducible divisors 
$D$ in $X$ such that $\num(D)=0$ is countable, and these divisors 
are rigid as well as their multiples. If $\alpha\in\cE$ is a 
pseudo-effective class of numerical dimension~$0$, then $\alpha$ 
is numerically equivalent to 
an effective $\bR$-divisor $D=\sum_{j\in J}\lambda_j D_j$,
for some finite subset $(D_j)_{j\in J}\subset\cD_0$ such that the 
cohomology classes $\{D_j\}$ are linearly independent and
some \hbox{$\lambda_j>0$}. If such a linear combination is of
numerical dimension~$0$, then so is any other linear combination of the
same divisors.
\endclaim

\proof. It is immediate from the definition that a pseudo-effective
class is of numerical dimension $0$ if and only if $\langle\alpha\rangle=0$,
in other words if $\alpha=N(\alpha)$. Thus $\alpha\equiv\sum\lambda_jD_j$
as described in 3.7, and since $\lambda_j\langle D_j\rangle\le 
\langle\alpha\rangle$, the divisors $D_j$ must themselves have numerical 
dimension $0$. There is at most one such divisor $D$ in any given cohomology
class in $NS(X)\cap\cE\subset H^2(X,\bZ)$, otherwise two such divisors 
$D\equiv D'$ would yield a blow-up $\mu:\wt X\to X$ resolving the 
intersection, and by taking
$\min(\mu^\star D,\mu^\star D')$ via $(3.4')$, we would find 
$\mu^\star D\equiv E+\beta$, $\beta\ne 0$, so that
$\{D\}$ would not be of numerical dimension $0$. This implies 
that there are at most countably many divisors of numerical 
dimension~$0$, and that these divisors are rigid as well as their
multiples.\qed
\medskip

The above general concept of numerical dimension leads to a very natural 
formulation of the abundance conjecture for non-minimal (K\"ahler) varieties.

\claim 3.8 Generalized abundance conjecture|For an arbitrary compact 
K\"ahler manifold~$X$, the Kodaira dimension should be equal to the
numerical dimension~:
$$\kappa(X)=\num(X):=\num(c_1(K_X)).$$
\endclaim

This appears to be a fairly strong statement. In fact, it is not difficult
to show that the generalized abundance conjecture would contain the
$C_{n,m}$ conjectures.

\claim 3.9 Remark|{\rm Using the Iitaka fibration, it is immediate to 
see that $\kappa(X)\leq\num(X)$.}
\endclaim

\claim 3.10 Remark|{\rm It is known that abundance holds in case
\hbox{$\num(X)=-\infty$} (if $K_X$ is not pseudo-effective, no
multiple of $K_X$ can have sections), or in case $\num(X)=n$. The
latter follows from the solution of the Grauert-Riemenschneider
conjecture in the form proven in [De85] (see also [DPa03]).

In the remaining cases, the most tractable situation is probably the
case when $\num(X)=0$.
In fact Theorem 3.7 then gives $K_X\equiv\sum\lambda_j D_j$ for some
effective divisor with numerically independent components, $\num(D_j)=0$. 
It follows that the $\lambda_j$ are rational and therefore
$$
K_X\sim \sum\lambda_j D_j+F\qquad
\hbox{where $\lambda_j\in\bQ^+$, $\num(D_j)=0$ and $F\in\Pic^0(X)$.}
\leqno(*)
$$
Especially, if we assume additionally that
$q(X)=h^{0,1}(X)$ is zero, then $mK_X$ is linearly equivalent
to an integral divisor for some multiple $m$, and it follows 
immediately that $\kappa(X)=0$. The case of a general projective
(or compact K\"ahler) manifold with \hbox{$\num(X)=0$} and positive 
irregularity $q(X)>0$ would be interesting to understand.}
\endclaim

The preceeding remarks at least give a proof up to dimension $4:$

\claim 3.11 Proposition|Let $X$ be a smooth projective $n$-fold with 
$n\leq 4$. If $\nu (X) = 0$, then $\kappa (X) = 0$. 
\endclaim

\proof. The proof is given in (9.1) in a slightly more general situation.
\qed

We will come back to abundance on 4-folds in sect. 9.

\section{\S4 The orthogonality estimate}

The goal of this section is to show that, in an appropriate sense,
approximate Zariski decompositions are almost orthogonal.

\claim 4.1 Theorem|Let $X$ be a projective manifold, and
let $\alpha=\{T\}\in\cE^\circ_\NS$ be a big class represented by
a K\"ahler current~$T$. Consider
an approximate Zariski decomposition
$$\mu_m^\star T_m = [E_m]+[D_m]$$
Then
$$
(D_m^{n-1}\cdot E_m)^2\le 20\,(C\omega)^n\big(\Vol(\alpha)-D_m^n\big)
$$
where $\omega=c_1(H)$ is a K\"ahler form and $C\ge 0$ is a constant such that
$\pm\alpha$ is dominated by $C\omega$ $($i.e., $C\omega\pm\alpha$ is nef$\,)$.
\endclaim

\proof. For every $t\in[0,1]$, we have
$$
\Vol(\alpha)=\Vol(E_m+D_m)\ge \Vol(tE_m+D_m).
$$
Now, by our choice of $C$, we can write $E_m$ as a difference of two
nef divisors
$$E_m=\mu^\star\alpha-D_m=\mu_m^\star(\alpha+C\omega)-(D_m+C\mu_m^\star
\omega).$$

\claim 4.2 Lemma|For all nef $\bR$-divisors $A$, $B$ we have
$$\Vol(A-B)\ge A^n - nA^{n-1}\cdot B$$
as soon as the right hand side is positive.
\endclaim

\proof. In case $A$ and $B$ are integral (Cartier) divisors, this is a 
consequence of the holomorphic Mores inequalities, [De01,8.4]. If $A$ and 
$B$ are $\bQ$-Cartier, we conclude by the homogeneity of the volume. The 
general case of $\bR$-divisors follows by approximation using the
upper semi-continuity of the volume [Bou02b, 3.1.26].\qed 

\claim 4.3 Remark|{\rm We hope that Lemma 4.2 also holds true on an
arbitrary K\"ahler mani\-fold for arbitrary nef (non neces\-sarily
integral) classes. This would follow from a generalization of
holomorphic Morse inequalities to non integral classes. However 
the proof of such a result seems technically much more involved than
in the case of integral classes.}
\endclaim
 
\claim 4.4 Lemma|Let $\beta_1,\ldots,\beta_n$ and 
$\beta'_1,\ldots,\beta'_n$ be nef classes on a compact K\"ahler manifold
$\wt X$ such that each difference $\beta'_j-\beta_j$ is pseudo-effective. Then
the $n$-th intersection products satisfy
$$
\beta_1\cdots\beta_n \le \beta'_1\cdots\beta'_n.
$$ 
\endclaim

\proof. We can proceed step by step and replace just one $\beta_j$ by
$\beta'_j\equiv\beta_j+T_j$ where $T_j$ is a closed positive
$(1,1)$-current and the other classes $\beta'_k=\beta_k$, $k\ne j$ are
limits of K\"ahler forms. The inequality is then obvious.\qed

\noindent
{\it End of proof of Theorem 4.1.}
In order to exploit the lower bound of the volume, we write
$$
tE_m+D_m=A-B,\qquad
A=D_m+t\mu_m^\star(\alpha+C\omega),\quad B=t(D_m+C\mu_m^\star\omega).
$$
By our choice of the constant $C$, both $A$ and $B$ are nef. Lemma 4.2 and 
the binomial formula imply
$$\eqalign{
\Vol(tE_m+D_m)&\ge A^n - nA^{n-1}\cdot B\cr
&=D_m^n+nt\,D_m^{n-1}\cdot\mu_m^\star(\alpha+C\omega)+\sum_{k=2}^n
t^k{n\choose k}D_m^{n-k}\cdot\mu_m^\star(\alpha+C\omega)^k\cr
&\phantom{{}=D_m^n}-nt\,D_m^{n-1}\cdot(D_m+C\mu_m^\star\omega)\cr
&\phantom{{}=D_m^n}-nt^2\sum_{k=1}^{n-1}t^{k-1}{n-1\choose k}
D_m^{n-1-k}\cdot\mu_m^\star(\alpha+C\omega)^k\cdot(D_m+C\mu_m^\star\omega).
\kern-4pt\cr}
$$
Now, we use the obvious inequalities
$$
D_m\le \mu_m^\star(C\omega),\qquad
\mu_m^\star(\alpha+C\omega)\le 2\mu_m^\star(C\omega),\qquad
D_m+C\mu_m^\star\omega\le 2\mu_m^\star(C\omega)
$$
in which all members are nef (and where the inequality${}\le$ means 
that the difference of classes is pseudo-effective). We use Lemma~4.4
to bound the last summation in the estimate of the volume, and in this
way we get
$$
\Vol(tE_m+D_m)\ge D_m^n+nt D_m^{n-1}\cdot E_m-nt^2\sum_{k=1}^{n-1}2^{k+1}
t^{k-1}{n-1\choose k}(C\omega)^n.
$$
We will always take $t$ smaller than $1/10n$ so that the last summation is 
bounded by $4(n-1)(1+1/5n)^{n-2}<4ne^{1/5}<5n$. This implies
$$
\Vol(tE_m+D_m)\ge D_m^n+nt\,D_m^{n-1}\cdot E_m-5n^2t^2(C\omega)^n.
$$
Now, the choice $t={1\over 10n}(D_m^{n-1}\cdot E_m)((C\omega)^n)^{-1}$ 
gives by substituting
$$
\eqalign{
{1\over 20}{(D_m^{n-1}\cdot E_m)^2\over (C\omega)^n}
&\le \Vol(E_m+D_m)-D_m^n\le \Vol(\alpha)-D_m^n\cr}
$$
(and we have indeed $t\le {1\over 10n}$ by Lemma 4.4), whence Theorem 4.1. 
Of course, the constant 20 is certainly not optimal.\qed

\claim 4.5 Corollary|If $\alpha\in\cE_\NS$, then the divisorial
Zariski decomposition $\alpha=N(\alpha)+\langle\alpha\rangle$
is such that
$$
\langle\alpha^{n-1}\rangle\cdot N(\alpha)=0.
$$
\endclaim

\proof. By replacing $\alpha$ by $\alpha+\delta c_1(H)$, one sees that
it is sufficient to consider the case where $\alpha$ is big. Then the 
orthogonality estimate implies
$$
(\mu_m)_\star(D_m^{n-1})\cdot (\mu_m)_\star E_m
= D_m^{n-1}\cdot(\mu_m)^\star (\mu_m)_\star E_m
\le D_m^{n-1}\cdot E_m \le C (\Vol(\alpha)-D_m^n)^{1/2}.
$$
Since $\langle\alpha^{n-1}\rangle = \lim (\mu_m)_\star(D_m^{n-1})$,
$N(\alpha)=\lim (\mu_m)_\star E_m$ and $\lim D_m^n=\Vol(\alpha)$,
we get the desired conclusion in the limit.\qed

\section{\S5 Proof of the duality theorem}

We want to prove that $\cE_\NS$ and $\SME(X)$ are dual (Theorem 2.2). 
By 1.4~(iii) we have in any case
$$\cE_\NS\subset (\SME(X))^\smallvee.$$ 
If the inclusion is strict, there is an element
$\alpha\in\partial\cE_\NS$ on the boundary of $\cE_\NS$ which is in
the interior of $\SME(X)^\smallvee$.

Let $\omega=c_1(H)$ be an ample class. Since $\alpha\in\partial\cE_\NS$, 
the class $\alpha+\delta\omega$ is big for every $\delta>0$, and since
$\alpha\in((\SME(X))^{\smallvee})^\circ$ we still have 
$\alpha-\varepsilon\omega\in(\SME(X))^{\smallvee}$ for
$\varepsilon>0$ small. Therefore
$$
\alpha\cdot \Gamma\ge\varepsilon\omega\cdot \Gamma\leqno(5.1)
$$
for every movable curve $\Gamma$. We are going to contradict (5.1).
Since $\alpha+\delta\omega$ is big, we have an approximate Zariski
decomposition
$$
\mu_\delta^\star(\alpha+\delta\omega)=E_\delta+D_\delta.
$$
We pick $\Gamma=(\mu_\delta)_\star(D_\delta^{n-1})$. By the Hovanskii-Teissier
concavity inequality
$$
\omega\cdot \Gamma\ge (\omega^n)^{1/n}(D_\delta^n)^{(n-1)/n}.
$$
On the other hand
$$
\eqalign{
\alpha\cdot \Gamma&=\alpha\cdot(\mu_\delta)_\star(D_\delta^{n-1})\cr
&=\mu_\delta^\star\alpha\cdot D_\delta^{n-1}\le
\mu_\delta^\star(\alpha+\delta\omega)\cdot D_\delta^{n-1}\cr
&=(E_\delta+D_\delta)\cdot D_\delta^{n-1}
=D_\delta^n+D_\delta^{n-1}\cdot E_\delta.\cr
}
$$
By the orthogonality estimate, we find
$$\eqalign{
{\alpha\cdot \Gamma\over
\omega\cdot \Gamma}
&\le {D_\delta^n + \big(20(C\omega)^n
(\Vol(\alpha+\delta\omega)-D_\delta^n)\big)^{1/2}\over
(\omega^n)^{1/n}(D_\delta^n)^{(n-1)/n}}\cr
&\le C'(D_\delta^n)^{1/n}+C''{(\Vol(\alpha+\delta\omega)-D_\delta^n)^{1/2}\over
(D_\delta^n)^{(n-1)/n}}.\cr}
$$
However, since $\alpha\in\partial\cE_\NS$, the class $\alpha$ cannot be big 
so
$$\lim_{\delta \to 0}D_\delta^n=\Vol(\alpha)=0.$$ 
We can also take $D_\delta$ to approximate $\Vol(\alpha+\delta\omega)$ 
in such a way that
$(\Vol(\alpha+\delta\omega)-D_\delta^n)^{1/2}$ tends to $0$ much faster than
$D_\delta^n$. Notice that \hbox{$D_\delta^n\ge\delta^n\omega^n$}, so in fact
it is enough to take
$$
\Vol(\alpha+\delta\omega)-D_\delta^n\le \delta^{2n}.
$$
This is the desired contradiction by (5.1).\qed
\bigskip

\claim 5.2 Remark|{\rm If holomorphic Morse inequalities were known also
in the K\"ahler case, we would infer by the same proof that
``$\alpha$ not pseudo-effective'' implies the existence of a blow-up
$\mu:\wt X\to X$ and a K\"ahler metric 
$\wt\omega$ on $\smash{\wt X}$ such that 
$\alpha\cdot\mu_\star(\wt\omega)^{n-1} < 0$. In the special case
when $\alpha=K_X$ is not pseudo-effective, we would expect the 
K\"ahler manifold $X$ to be covered by rational curves.
The main trouble is that characteristic $p$ techniques are no longer
available. On the other hand it is tempting to approach the
question via techniques of symplectic geometry~: } \endclaim

\claim 5.3 Question|Let $(M,\omega)$ be a compact real symplectic manifold.
Fix an almost\break complex structure $J$ compatible with $\omega$, and
for this structure, assume that\break \hbox{$c_1(M)\cdot\omega^{n-1}>0$}. 
Does it follow that $M$ is covered by rational $J$-pseudoholomorphic curves~?
\endclaim

\section{\S6 Non nef loci}

Following [Bou02b], we introduce the concept of non-nef locus of an arbitrary
pseudo-effective class. The details differ a little bit here
(and are substantially simpler) because the scope is limited to
compact K\"ahler manifolds.

\claim 6.1 Definition|Let $X$ be a compact K\"ahler manifold,
$\omega$ a K\"ahler metric, and $\alpha\in\cE$ a pseudo-effective class. 
We define the non-nef locus of $\alpha$ to be
$$
L_\nnef(\alpha)=\bigcup_{\delta>0} \bigcap_T \mu(|E|)
$$
for all log resolutions $\mu^\star T=[E]+\beta$ of positive currents
$T\in\{\alpha+\delta\omega\}$ with logarithmic singularities,
$\mu:\smash{\wt X}\to X$, and $\mu(|E|)$ is the set-theoretic image
of the support of $E$.
\endclaim

It should be noticed that the union in the above definition can be restricted
to any sequence $\delta_k$ converging to $0$, hence $L_\nnef(\alpha)$ is
either an analytic set or a countable union of analytic sets. The results
of [De92] and [Bou02b] show that 
$$
L_\nnef(\alpha)=\bigcup_{\delta>0} \bigcap_T E_+(T)\leqno(6.1')
$$
where $T$ runs over the set $\alpha[-\delta\omega]$ of all $d$-closed
real $(1,1)$-currents $T\in\alpha$ such that $T\ge-\delta\omega$, and
$E_+(T)$ denotes the locus where the Lelong numbers of $T$ are
strictly positive. The latter definition $(6.1')$ works even in the
non K\"ahler case, taking $\omega$ an arbitrary positive hermitian
form on~$X$. By [Bou02b], there is always a current $T_{\min}$ which
achieves minimum singularities and minimum Lelong numbers among
all members of $\alpha[-\delta\omega]$, hence
$\bigcap_T E_+(T)=E_+(T_{\min})$.

\claim 6.2 Theorem|Let $\alpha\in\cE$ be a pseudo-effective class.
Then $L_\nnef(\alpha)$ contains the union of all irreducible algebraic 
curves $C$ such that $\alpha \cdot C<0$.
\endclaim

\proof. If $C$ is an irreducible curve not contained in 
$L_\nnef(\alpha)$, the definition implies that for every $\delta>0$
we can choose a positive current $T\in\{\alpha+\delta\omega\}$ and
a log-resolution $\mu^\star T=[E]+\beta$ such that $C\not\subset\mu(|E|)$.
Let $\wt C$ be the strict transform of $C$ in $\wt X$, so that $C=\mu_\star
\wt C$. We then find
$$
(\alpha+\delta\omega)\cdot C=([E]+\beta)\cdot\wt C\ge 0
$$
since $\beta\ge 0$ and $\wt C\not\subset |E|$. This is true for
all $\delta>0$ and the claim follows.
\medskip

\claim 6.3 Remark|{\rm One may wonder, at least when $X$ is projective
and $\alpha\in\cE_\NS$, whether $L_\nnef(\alpha)$ is actually equal
to the union of curves $C$ such that $L\cdot C<0$ (or the ``countable
Zariski closure'' of such a union). Unfortunately, this is not true, 
even on surfaces. The following simple example
was shown to us by E.~Viehweg. Let $Y$ be a complex algebraic surface 
possessing a big line bundle $F$ with a curve $C$ such that $F\cdot C<0$ 
as its base locus (e.g.\ $F=\pi^\star\cO(1)+E$ for the blow-up $\pi:Y\to\bP^2$
of $\bP^2$ in one point, and $C=E={}$exceptional divisor). Then 
take finitely many points $p_j\in C$, $1\le j\le N$, and blow-up these 
points to get a modification $\mu:X\to Y$. We select
$$
L=\mu^\star F+\wh C+2\sum E_j=\mu^\star(F+C)+\sum E_j
$$
where $\wh C$ is the strict transform of $C$ and $E_j=\mu^{-1}(p_j)$.
It is clear that the non nef locus of $\alpha=c_1(L)$ must be equal
to $\smash{\wh C}\cup\bigcup E_j$, although 
$$
L\cdot\wh C=(F+C)\cdot C + N>0
$$
for $N$ large. This example shows that the set of $\alpha$-negative curves 
is not the appropriate tool to understand the non nef locus.}
\endclaim

\section{\S7 Pseudo-effective vector bundles}

In this section we consider pseudo-effective and almost nef vector
bundles as introduced in [DPS00]. As an application, we obtain interesting 
informations concerning the tangent bundle of Calabi-Yau manifolds. 
First we recall the relevant definitions.

\claim 7.1 Definition|Let $X$ be a compact K\"ahler manifold and $E$ a
holomorphic vector bundle on $X$. Then $E$ is said to be {\it
pseudo-effective} if the line bundle $\cO_{\bP(E)}(1)$ is 
pseudo-effective on the projectivized bundle $\bP(E)$ of hyperplanes
of $E$, and if the projection $\pi(L_\nnef(\cO_{\bP(E)}(1)))$
of the non-nef locus of $\cO_{\bP(E)}(1)$ onto $X$ does not cover all 
of~$X$.
\endclaim

This definition would even make sense on a general compact complex manifold,
using the general definition of the non-nef locus in [Bou02b]. On the other
hand, the following proposition gives an algebraic characterization of
pseudo-effective vector bundles in the projective case.

\claim 7.2 Proposition|Let $X$ be a projective manifold. A holomorphic
vector bundle $E$ on $X$ is {\it pseudo-effective} if and only if
for any given ample line bundle $A$ on $X$ and any positive integers
$m_0$, $p_0$, the vector bundle
$$
S^p((S^mE)\otimes A)
$$
is generically generated $($i.e.\ generated by its global sections
on the complement $X\ssm Z_{m,p}$ of some algebraic set $Z_{m,p}\ne X)$
for some $[$resp.\ every$]$ $m \geq m_0 $ and $p \geq p_0$.
\endclaim

\proof. If global sections as in the statement of 7.2
exist, they can be used to define a singular hermitian metric $h_{m,p}$ on 
$\cO_{\bP(E)}(1)$ which has poles contained in $\pi^{-1}(Z_{m,p})$ and whose
curvature form satisfies $\Theta_{h_{m,p}}(\cO_{\bP(E)}(1))\ge 
-{1\over m}\pi^* \Theta(A)$. Hence, by selecting suitable
integers $m=M(m_0,p_0)$ and $p=P(m_0,p_0)$, we find that 
$\cO_{\bP(E)}(1)$ is pseudo-effective (its first Chern class is 
a limit of pseudo-effective classes), and that
$$
\pi(L_\nnef(\cO_{\bP(E)}(1)))\subset \bigcup_{m_0}\bigcap_{p_0}
Z_{m,p}\subsetneq X.
$$
Conversely, assume that $\cO_{\bP(E)}(1)$ is pseudo-effective and admits
singular hermitian metrics $h_\delta$ such that 
$\Theta_{h_\delta}(\cO_{\bP(E)}(1))\ge-\delta\wt\omega$ and
$\pi(\Sing(h_\delta))\subset Z_\delta\subsetneq X$ (for some K\"ahler
metric $\wt\omega$ on $\bP(E)$ and arbitrary small $\delta>0$). 
We can actually take $\omega=\Theta(A)$ and
$\wt\omega=\varepsilon_0\Theta_{h_0}(\cO_{\bP(E)}(1))+\pi^*\omega$
with a given smooth hermitian metric $h_0$ on $E$ and
$\varepsilon_0\ll 1$.  An easy calculation shows that the linear
combination $h'_\delta=
\smash{h_\delta^{1/(1+\delta\varepsilon_0)}h_0^{\delta\varepsilon_0}}$ 
yields a metric on $\cO_{\bP(E)}(1)$ such that
$$
\Theta_{h'_\delta}(\cO_{\bP(E)}(1))\ge -\delta\pi^*\Theta(A).
$$
By taking $\delta={1/2m}$ and multiplying by $m$, we find
$$
\Theta(\cO_{\bP(E)}(m)\otimes \pi^*A) \ge {1\over 2}\pi^*\Theta(A)
$$
for some metric on $\cO_{\bP(E)}(m)\otimes \pi^*A$ which is smooth
over $\pi^{-1}(X\ssm Z_\delta)$. The standard theory of $L^2$ estimates
for bundle-valued $\overline\partial$-operators can be used to produce the 
required sections, after we multiply $\Theta(A)$ by a sufficiently
large integer $p$ to compensate the curvature of $-K_X$. The sections
possibly still have to vanish along the poles of the metric, but
they are unrestricted on fibers of $\bP(S^mE)\to X$ which do not meet
the singularities.\qed
\medskip

Note that if $E$ is pseudo-effective, then $\cO_{\bP(E)}(1)$ is
pseudo-effective and $E$ is {\it almost nef} in the following sense
which is just the straightforward generalization from the line bundle
case.

\claim 7.3 Definition|Let $X$ be a projective manifold and $E$ a
vector bundle on~$X$. Then $E$ is said to be {\it almost nef}, if
there is a countable family $A_i $ of proper subvarieties of $X$ such
that $E \vert C$ is nef for all $C \not \subset \bigcup _{i} A_i$. 
Alternatively, $E$ is almost nef if there is no covering family of 
curves such that $E$ is non-nef on the general member of the family. 
\endclaim

Observe that  $E$ is almost nef if and only if
$\cO_{\bP(E)}(1)$ is almost nef and $\cO_{\bP(E)}(1)$ is nef on the 
general member of any family of curves
in $\bP(E)$ whose images cover $X$. 
Hence Theorem 2.2 yields

\claim 7.4 Corollary| Let $X$ be a projective manifold and $E$ a
holomorphic vector bundle on $X$.  If $E$ is almost nef, 
then $\cO_{\bP(E)}(1)$ is pseudo-effective. Thus 
for some $[$or any$]$ ample line bundle $A$, there are positive 
numbers $m_0$ and $p_0$
such that
$$
H^0(X,S^p((S^mE) \otimes A)) \ne 0
$$
for all $m \geq m_0 $ and $p \geq p_0$.
\endclaim

One should notice that it makes a big difference to assert just the
existence of a non zero section, and to assert the existence of
sufficiently many sections guaranteeing that the fibers are
generically generated.  It is therefore natural to raise the following
question.

\claim 7.5 Question| Let $X$ be a projective manifold and $E$ a vector
bundle on $X$. Suppose that $E$ is almost nef. Is $E$ always 
pseudo-effective in the sense of Definition~$7.1$~?
\endclaim

This was stated as a theorem in [DPS01, 6.3], but the proof given
there was incomplete.  The result now appears quite doubtful to us.
However, we give below a positive answer to Question~7.5 in case of a
rank $2-$bundle $E$ with $c_1(E) = 0$ (conjectured in [DPS01]), and
then apply it to the study of tangent bundles of K3-surfaces.

\claim 7.6 Theorem|Let $E$ be an almost nef vector bundle of rank at most 
$3$ on a projective manifold~$X$. Suppose that $\det E \equiv 0$. 
Then $E$ is numerically flat. 
\endclaim

\proof. Recall (cf.\ [DPS94]) that a vector bundle $E$ is said to be
numerically flat if it is nef as well as its dual (or, equivalently,
if $E$ is nef and $\det E$ numerically trivial); also, $E$
is numerically flat if and only if $E$ admits a filtration by
subbundles such that the graded pieces are unitary flat vector
bundles.  By [Ko87, p.115], $E$ is unitary flat as soon as
$E$ is stable for some polarization and $c_1(E)=c_2(E)=0$.

Under our assumptions, $E$ is necessarily semi-stable since
semi-stability with respect to a polarization $H$ can be tested
against a generic complete intersection curve, and we know that $E$ is
nef, hence numerically flat, on such a curve. Therefore (see also [DPS01,
6.8]) we can assume without loss of generality that $\dim X=2$ and
that $E$ is stable with respect to all polarizations, and it is enough
to show in that case that $c_2(E)=0$. Since $E$ is almost nef, $E$ is
nef, hence numerically flat, on all curves except for at most a
countable number of curves, say $(\Gamma_j)_{j\in\bN}$.
\\
First suppose that $E$ has rank $2.$
Then the line bundle $\cO(1)$ on $\bP(E) $ is immediately seen to
be nef on all but a countable number of curves. In fact, the only
curves on which $\cO(1)$ is negative are the sections over the curves
$\Gamma_j$ with negative self-intersection in $\bP(E \vert \Gamma_j)$.
Now take a general hyperplane section $H$ on $\bP(E).$ Then $H$ does
not contain any of these bad curves and therefore $\cO(1)$ is nef on
$H.$ Hence
$$
c_1(\cO(1))^2 \cdot H \geq 0.
$$ 
Now - up to a multiple - $H$ is of the form $H = \cO(1) \otimes \pi^*(G)$ 
so that 
$$
c_1(\cO(1))^3 + c_1(\cO(1))^2 \cdot \pi^*(G) \geq 0.
$$
Since $c_1(\cO(1))^3 = c_1(E)^2 - c_2(E) = - c_2(E)$ and $c_1(\cO(1))^2 
\cdot \pi^*(G) = c_1(E) \cdot G = 0,$ we conclude $c_2(E) = 0.$
\\
If $E$ has rank $3,$ we need to argue more carefully, because now
$\cO(1)$ is non-nef on the surfaces $S_j = \bP(E \vert C_j)$ so that
$\cO(1) $ might be non-nef on a general hyperplane section $H.$ We
will however show that this can be avoided by choosing carefully the
linear system $\vert H \vert.$ To be more precise we fix $G$ ample on
$X$ and look for
$$
H \in \vert \cO(1) + \pi^*(mG) \vert
$$
with $m \gg 0, $ so that $\cO(1)$ is nef on $H \cap S_j$ for all
$j.$ Given that $\cO(1) \vert H$ and we can argue as in the previous
case to obtain $c_2(E) = 0.$ Of course for a general choice of $H$,
all curves $H \cap S_j$ will be irreducible (but possibly singular
since $C_j$ might be singular).  Now fix $j$ and set $\tilde C = H
\cap S_j, $ a section over $C = C_j.$ Let $V \subset E_C$ be the
maximal ample subsheaf (see [PS02]). Then we obtain a vector bundle
sequence
$$
0 \to V \to E_C \to F \to 0
$$
and we may assume that $F$ has rank $2,$ because otherwise $\cO(1)$
is not nef only on one curve over $C.$ Now $\tilde C$ induces an exact
sequence
$$
0 \to \cO_C(-mG) \to F \to F' \to 0
$$
and therefore $\cO(1) \vert \tilde C$ is nef iff $c_1(F') \geq 0.$ 
This translates into $c_1(F) + m(G \cdot C) \geq 0.$ 
Now let $t_0$ be the nef value of $E$ with respect to $G$, i.e.\ $E(t_0G)$ 
is nef but not ample. Then $F(t_0G)$ is nef, too, so that
$c_1(F) \geq -2t_0 (G \cdot C).$ In total 
$$
c_1(F') \geq (m-2t_0) (G \cdot C),
$$
hence we choose $m \geq 2t_0$ and for this choice $\cO(1) \vert H$ is nef. 
\qed
  
As a corollary we obtain

\claim 7.7 Theorem|Let $X$ be a projective 
K3-surface or a Calabi-Yau 3-fold. Then the tangent bundle $T_X$ is not 
almost nef, and there exists a covering family $(C_t)$ of (generically
irreducible) curves such that $T_X \vert C_t$ is not nef for general $t$.
\endclaim

In other words, if $c_1(X) = 0$ and $T_X$ is almost nef, then a finite
\'etale cover of $X$ is abelian. 
One should compare this with Miyaoka's theorem that $T_X \vert C$ is
nef for a smooth curve $C$ cut out by hyperplane sections of
sufficiently large degree. Note also that $T_X \vert C$ being not nef
is equivalent to say that $T_X \vert C$ is not semi-stable.
We expect that (7.7) holds in general for Calabi-Yau manifolds of any 
dimension.

\proof. Assume that $T_X$ is almost nef. Then by 7.6, $T_X$ is numerically 
flat. In particular $c_2(X) = 0$ and hence $X$ is an \'etale quotient
of a torus.\qed

\noindent We will now improve (7.7) for K3-surfaces; namely if $X$ is 
a projective K3-surface, then already $\cO_{\bP(T_X)}(1)$ should be
non-pseudo-effective.  In other words,  let $A$ be a fixed ample
divisor on $X$. Then for all positive integers $m$ there exists a
positive integer $p$ such that
$$
H^0(X,S^p((S^mT_X) \otimes A)) = 0.
$$
This has been verified in [DPS00] for the general quartic in
$\bP_3$ and below for any K3-surface.

\claim 7.8 Theorem| Let $X$ be a projective K3-surface and 
$L = \cO_{\bP(T_X)}(1)$. Then $L$ is not pseudo-effective.
\endclaim

\proof.  Suppose that $L$ is pseudo-effective and consider the divisorial 
Zariski decomposition ([Bou02b], cf.\ also 3.5~(iv))
$$
L = N + Z
$$
with $N$ an effective $\bR$-divisor and $Z$ nef in codimension $1$.
Write $N = aL + \pi^*(N') $ and $Z = bL + \pi^*(Z'). $ Let $H$ be very
ample on $S$. By restricting to a general curve $C$ in $\vert H_H
\vert $ and observing that $T_X \vert C$ is numerically flat,
we see that $L \vert \pi^{-1}(C)$ is nef (but not ample), hence
$$ N' \cdot C = 0.$$
Thus $N' = 0$ since $H$ is arbitrary. If $a > 0,$ then some $mL$ would be
effective, i.e. $S^mT_X$ would have a section, which is known not to be
the case. Hence $a = 0$ and $L$ is nef in codimension 1
so that $L$ can be negative only on finitely many curves. This
contradicts 7.7.  \qed

\section{\S8 Partial nef reduction} 

In this section we construct reduction maps for pseudo-effective line
bundle which are zero on large families of curves. This will be
applied in the next section in connection with the abundance problem.

\claim 8.1 Notation| {\rm Let $(C_t)_{t \in T}$ be a covering family
of (generically irreducible) curves (in particular $T$ is
irreducible and compact).  Then $(C_t)$ is said to be a {\it
connecting family} if and only if two general $x,y$ can be joined
by a chain of $C_t$. We also say that $X$ is $(C_t)$-connected.}
\endclaim

Using Campana's reduction theory [Ca81,94], we obtain immediately

\claim 8.2 Theorem| Let $X$ be a projective manifold and $L$ a
pseudo-effective line bundle on $X$. Let $(C_t) $ be a covering family
with $L \cdot C_t = 0$. Then there exists an almost holomorphic
surjective meromorphic map $f: X \merto Y$ such that the general
(compact) fiber of $f$ is $(C_t)$-connected. $f$ is called the partial
nef reduction of $L$ with respect to $(C_t)$.
\endclaim

\claim 8.3 Definition| Let $L$ be a pseudo-effective line bundle on
$X$. The minimal number which can be realised as $\dim Y$ with a
partial nef reduction $f: X \merto Y$ with respect to $L$ is denoted
$p(L)$.  If there is no covering family $(C_t)$ with $L \cdot C_t = 0$, 
then we set $p(L) = \dim X$.
\endclaim

\claim 8.4 Remark| {\rm $p(L) = 0$ if and only if there exists a connecting 
family $(C_t)$ such that $L \cdot C_t = 0$. If moreover $L$ is nef, then
$p(L) = 0$ if and only if $L \equiv 0$ [Workshop].}
\endclaim

\claim 8.5 Proposition| Let $X$ be a projective manifold, $L$ a
pseudo-effective line bundle and $(C_t)$ a connecting family. If $L
\cdot C_t = 0$, then $\kappa(L) \leq 0$.  
\endclaim

\proof.  Supposing the contrary we may assume that $h^0(L) \geq 2$.
Let $C = C_t$ be a general member of our family. Then we find a
non-zero $s \in H^0(L) $ such that $s \vert C = 0$. If $A \subset X$
is any algebraic set, we let $G(A) $ be the union of all $y$ which can
be joined with $A$ by a single $C_t$. Now the general $C_t$ through a
general point must be irreducible, hence $G(C)$ is generically filled
up by irreducible $C_t$, and we conclude that $s \vert G(C) = 0$.
Define $G^k(C) = G(G^{k-1}(C))$. Then by induction we obtain $s \vert
G^k(C) = 0$ for all $k$. Since $(C_t)$ is connecting, we have
$G^{\infty} = X$, hence $s = 0$. Therefore $h^0(L) \leq 1$.
\qed

It is also interesting to look at covering families $(C_t)$ of {\it ample} 
curves. Here ``ample'' means that the dual of the conormal
sheaf modulo torsion is ample (say on the normalization). Then we have
the same result as in (8.5) which is prepared by

\claim 8.6 Lemma| Let $X$ be a projective manifold, $C \subset X$ an
irreducible curve with normalization $f: \tilde C \to C$ and ideal
sheaf $\cI$.  Let $L$ be a line bundle on $X$. Then there exists a
positive number $c$ such that for all $t \geq 0$:
$$
h^0(X,L^t) \leq \sum_{k=0}^{ct}h^0(f^*(S^k(\cI/\cI^2/{\rm tor})
\otimes L^t)).
$$
\endclaim

\proof.  Easy adaptation of the proof of (2.1) in [PSS99].
\qed

\claim 8.7 Corollary| Let $X$ be a projective manifold and $C \subset
X$ be an irreducible curve with normalisation $f: \tilde C \to C$
such that $f^*(\cI /\cI^2)^*$ is ample. Let $L$ a line bundle with $L
\cdot C_t = 0$.  Then $\kappa(L) \leq 0$. In particular this holds
for the general member of an ample covering family.
\endclaim

\proof.  By (8.5) it suffices to show that
$$
h^0(f^*(S^k(\cI/\cI^2/{\rm tor}) \otimes L^t)) = 0
$$
for all $k \geq 1$. This is however clear since by assumption
$f^*(\cI/\cI^2)^*$ is an ample bundle.  \qed

\claim 8.8 Corollary| Let $X$ be a smooth projective threefold with
$K_X$ pseudo-effective. If there is a ample covering family or a
connecting family $(C_t)$ such that $K_X \cdot C_t = 0$, then 
$\kappa(X) = 0$.
\endclaim

\proof.  By (8.5) resp.\ (8.7) we have $\kappa(X) \leq 0$. Suppose
that $\kappa(X) = -\infty$. Then $X$ is uniruled by Miyaoka's
theorem.  Thus $K_X$ is not pseudo-effective.
\qed

Although we will not need it later, we will construct a nef reduction
for pseudo-effective line bundles, generalizing to a certain extent
the result of [Workshop] for nef line bundle (however the result is
weaker). A different type of reduction was constructed in
[Ts00],[Ec02].

\claim 8.9 Theorem| Let $L$ be a pseudo-effective line bundle on a
projective manifold $X$. Then there exists an almost holomorphic
meromorphic map $f: X \merto Y$ such that
\smallskip 
\item {\rm (i)} general points on the general fiber of $F$ can be
  connected by a chain of $L$-trivial irreducible curves.
\smallskip
\item {\rm (ii)} if $x \in X$ is general and $C$ is an irreducible
  curve through $x$ with $\dim f(C) > 0$, then $L \cdot C > 0$.
\vskip0pt
\endclaim

\proof. The argument is rather standard: start with a covering
$L$-trivial family $(C_t)$ and build the nef partial quotient 
$h: X \merto Z$ (if the family does not exist, put $f = id)$. Now take
another covering $L$-trivial family $(B_s)$ (if this does not exist,
just stop) with partial nef reduction $g$. For general $z \in Z$ let
$F_z$ be the set of all $x \in X$ which can be joined with the fiber
$X_z$ by a chain of curves $B_s$.  In other words, $F_z$ is the
closure of $g^{-1}(g(X_z))$. Now the $F_z$ define a covering family
(of higher-dimensional subvarieties) which defines by Campana's
theorem a new reduction map. After finitely many steps we arrive at
the map we are looking for.  \qed

Finally we show that covering families which are interior points in
the movable cone are connecting:

\claim 8.10 Theorem| Let $X$ be a projective manifold and $(C_t)$ a
covering family. Suppose that $[C_t]$ is an interior point of the
movable cone $\cM$. Then $(C_t)$ is connecting.
\endclaim

\proof.  Let $f: X \merto Z$ be the reduction of the family $(C_t)$.
If the family is not connecting, then $\dim Z > 0$.  Let $\pi: \tilde
X \to X$ be a modification such that the induced map $\tilde f:
\tilde X \to Z$ is holomorphic. Let $A$ be very ample on $Z$ and put
$L = \pi_*(\tilde f^*(A))^{**}$. Then $L$ is an effective line bundle
on $X$ with $L \cdot C_t = 0$ since $L$ is trivial on the general
fiber of $f$, this map being almost holomorphic. Hence $[C_t]$ must be
on the boundary of $\cM$.
\qed

The converse of (8.10) is of course false: consider the family of
lines $l$ in $\bP_2$ and let $X$ be the blow-up of some point in
$\bP_2$.  Let $(C_t)$ be the closure of the family of preimages of
general lines. This is a connecting family, but if $E$ is the
exceptional divior, then $E \cdot C_t = 0$. So $(C_t)$ cannot be in
the interior of $\cM$.

\section{\S9 Towards abundance}

In this section we prove that a smooth projective 4-fold $X$ with
$K_X$ pseudo-effective and with the additional property that $K_X$ is
$0$ on some covering family of curves, has $\kappa(X) \geq 0$. In the
remaining case that $K_X$ is positive on all covering families of
curves one expects that $K_X$ is big.

\claim 9.1 Proposition| Let $X$ be a smooth projective 4-fold with
$K_X$ pseudo-effective. Suppose that there exists a dominant rational
map $f: X \merto Y$ to a projective manifold $Y$ with $\kappa(Y)\geq 0$ 
$($and $ 0 < \dim Y < 4)$.  Then $\kappa(X) \geq 0$.
\endclaim

\proof.  We may assume $f$ holomorphic with general fiber $F$.  If
$\kappa(F) = - \infty$, then $F$ would be uniruled, hence $X$ would
be uniruled. Hence $\kappa(F) \geq 0$.  Now $C_{n,n-3}$, $C_{n,n-2}$
and $C_{n,n-1}$ hold true, see e.g.\ [Mo87] for further references.
This gives
$$
\kappa(X) \geq \kappa(F) + \kappa(Y) \geq 0.
$$ 
\qed

\claim 9.2 Corollary| Let $X$ be a smooth projective 4-fold with $K_X
$ pseudo-effective. Let $f: X \merto Y$ be a dominant
rational map $(0 < \dim Y < 4)$ with $Y$ not rationally connected.
Then $\kappa(X) \geq 0$. \endclaim

\proof.  If $\dim Y \leq 2$, this is immediate from (9.1). So let
$\dim Y = 3$. Since we may assume $\kappa(Y) = - \infty$, the
threefold $Y$ is uniruled. Let $h: Y \merto Z$ be the
rational quotient; we may assume that $h$ is holomorphic and $Z$
smooth. Since $Y$ is not rationally connected, $\dim Z \geq 1$.  Then
$q(Z) \geq 1$, otherwise $Z$ would be rational and hence $Y$
rationally connected by Colliot-Th\'el\`ene [CT86], see also
Graber-Harris-Starr [GHS03], and we conclude by (9.1).
\qed

\claim 9.3 Conclusion| In order to prove $\kappa(X) \geq 0 $ in case
of a dominant rational map $f: X_4 \merto Y$, we may assume
that $Y$ is a rational curve, a rational surface or a rationally
connected 3-fold.
\endclaim

\claim 9.4 Proposition| Let $X$ be a smooth projective 4-fold with
$K_X$ pseudo-effective. If $p(K_X) = 1$, then $\kappa(X) \geq 0$.
\endclaim

\proof.  By assumption we have a covering family $(C_t)$ with $K_X
\cdot C_t = 0$, whose partial nef reduction is a holomorphic map $f: X
\to Y$ to a curve $Y$.  By (9.3) we may assume $Y = \bP_1$.  We
already saw that $\kappa(F) \geq 0$, however by (8.5) we even have
$\kappa(F) = 0$. Choose $m$ such that $h^0(mK_F) \ne 0$ for the
general fiber $F$ of $f$. Thus $f_*(mK_X)$ is a line bundle on $Y$,
and we can write
$$
mK_X = f^*(A) + \sum a_i F_i + E \leqno (*)
$$
where $F_i$ are fiber components and $E$ surjects onto $Y$ with
$h^0(\cO_X(E)) = 1$. The divisor $E$ comes from the fact that $F$ is
not necessarily minimal; actually $E \vert F = mK_F$. We will get rid
of $E$ by the following construction. Let $Y_0 \subset Y$ be the
largest open set such that $f$ is smooth over $Y_0;$ let $X_0 =
f^{-1}(Y_0)$. By [KM92] we have a birational model $f_0': X_0' \merto Y_0$ 
via a sequence of relative contractions and relative flips such that the
fibers $F'$ of $f_0'$ are minimal, hence $nK_{F'} = \cO_{F'}$ for
suitable $n$. Now compactify. Thus we may assume that the general
fiber of $f$ is minimal, paying the price that $X$ might have terminal
singularities. However these singularities don't play any role since below 
we will argue on a general surface in $X$, which automatically does not meet
the singular locus of $X$, this set being of dimension at most~$1$.
\bigskip

\noindent
In particular we have $E = 0$ in $(*)$. By enlarging $m$ we may also
assume that the support of $\sum a_i F_i$ does not contain any fiber
and also that $mK_X$ is Cartier.  Now let $S \subset X$ be a surface
cut out by 2 general hyperplane sections. Let $L = mK_X \vert S$.
Denoting $G_i = F_i \vert S$ and $g = f \vert S$, we obtain from $(*)$
$$
L = g^*(A) + \sum a_i G_i. \leqno (**)
$$
On the other hand, we consider the divisorial Zariski decomposition $L
\equiv N + Z$ constructed in [Bou02b]; see (3.5). Here $N$ is an effective
$\bR$-divisor covering the non-nef locus of $L$ and $Z$ is an
$\bR$-divisor which is nef in codimension 1. Let $l$ be a general fiber of $g$.
Then $L \cdot l = 0$
and thus
$$
N \cdot l = Z \cdot l = 0.
$$
So $N$ is contained in fibers of $g$ and $ Z = f^*(\cO_Y(a))$,
[Workshop,2.11]; moreover $a \geq 0$.  Comparing with $(**)$ we get $Z
= f^*(A)$ and thus $A$ is nef. Hence $(*)$ gives $\kappa(X) \geq 0$.
\qed

\claim 9.5 Remark| {\rm Proposition 9.4 also holds in dimension 5. In
fact, in view of $(C_{n,1})$, only two things to be observed. The
first is that $\kappa (F) \geq 0$. But this follows from (9.9)
below. The second is that in case of 4-dimensional fibers we cannot
apply [KM92]. However it is not really necessary to use [KM92]. We
can also argue as follows. Let $E' = E \vert S.$ Then necessarily
$E' \subset N$ so that $g^*(A) + \sum a_i G_i = N' + Z$ is
pseudo-effective, too, and we conclude as in (9.4).  \smallskip
\noindent The same remark also applies to the next proposition
(9.6). }
\endclaim

\claim 9.6 Proposition| Let $X$ be a smooth projective 4-fold with
$K_X $ pseudo-effective. If $p(K_X) = 2$, then $\kappa(X) \geq 0$.
\endclaim

\proof. By (9.3) we may assume that we have a holomorphic partial nef
reduction $f: X \to Y = \bP_2$. Again $\kappa(F) = 0$ for the
general fiber $F$. Then for a suitable large $m$ we have a
decomposition
$$
mK_X = f^*(A) + \sum a_i F_i + E_1-E_2 + D  \leqno (*)
$$
where the support of $\sum a_i F_i$ does not contain the support of
any divisor of the form $f^{-1}(C)$, but $f(F_i)$ is 1-dimensional for
all $i$, where $E_i$ are effective with $\dim f(E_i) = 0$, and where
$D$ is effective, projecting onto $Y$ with $D \vert F = mK_F$ for the
general fiber.  The divisor $E_2 $ arises from the fact that
$f_*(kK_X) $ might not be locally free, but only torsion free.
\smallskip

\noindent Writing $A = \cO_Y(a)$, we are going to prove that $a \geq 0$.  
Let $l \subset  Y$ be a general line and $X_l = f^{-1}(l)$. Let 
$g = f \vert X_l$, $A' = A \vert l$, $F'_i = F_i \cdot X_l$,
$D' = D \cdot X_l$ and $L = mK_X \vert X_l$. Then $(*)$ gives
$$
L = g^*(A') + \sum a_i F_i' + D'. 
$$
Passing to a suitable model for $X_l \to l$ as in the proof of
(9.4), we may assume that $D' = 0$.  Now consider Boucksom's
divisorial Zariski decomposition $L = N + Z $ as in the proof of
(9.4). Then we conclude as before that $Z = g^*(A') $ and that $N =
\sum a_i F_i'$. Hence $A'$ is nef and thus $a \geq 0$ (this can also
be verified easily without using the divisorial Zariski
decomposition).
\bigskip 

\noindent
Going back to $(*)$, the only remaining difficulty is the presence of
the negative summand $E_2$. This requires the following
considerations.  We write
$$
f_*(mK_X) = \cI_Z \otimes \cO(a)
$$
with a finite set $Z;$ being defined by $\cI_Z = f_*(\cO_X(E_2)).$
\smallskip \noindent In a first step we claim that $\cI_Z \otimes
\cO(a)$ is pseudo-effective in the sense that $ N(\cI_Z^k \otimes
\cO(ka+1)) $ has sections for large $k$ and $N.$ Let $\sigma: \hat Y
\to Y$ be a birational map with $\hat Y$ smooth and $\hat X$ the
normalization of the fiber product $X \times_Y \hat Y$ with induces
maps $\tau: \hat X \la X$ and $\hat f: \hat X \to \hat Y$ such that
$\hat f $ is flat, in particular equidimensional.  Adopting the
arguments from above, $\hat f_* \tau^*(mK_X)$ is a pseudo-effective
line bundle. Now
$$ \sigma_* \hat f_* \tau^*(mK_X) = f_* \tau_* \tau^*(mK_X) = f_*(mK_X),$$
hence $\cI_Z \otimes \cO(a)$ is clearly pseudo-effective. 
\smallskip

\noindent
Let $x \in Y = \bP_2$ be general and let $\rho: \tilde Y \to Y$ be the
blow-up of $x,$ inducing a $\bP_1-$ bundle $\pi: \tilde Y \to \bP_1.$
Let $F$ be a fiber of $\pi,$ i.e. a line in $\bP_2.$ Then the
pseudo-effectivity of $\cI_Z \otimes \cO(a) $ yields
$$ H^1(F,\cI_Z \otimes \cO(a) \vert F) = 0$$
from which we get the vanishing $R^1\pi_*(\cI_Z \otimes \rho^*(\cO(a)) = 0.$ 
Thus we have an exact sequence
$$ 0 \la \pi_*(\cI_Z \otimes \rho^*(\cO(a))) \la \pi_*(\rho^*(\cO(a))) 
\la \cO_R \la 0 $$
where $R \simeq Z.$ 
\smallskip

\noindent
Now suppose $H^0(\cI_Z \otimes \cO(a)) = 0$ (otherwise we are done).
Since $\pi_*(\rho^*(\cO(a)) = S^a(\cO \oplus \cO(1)), $ we deduce by
taking $H^0$ that
$$ l(Z) = {{a(a+1)} \over {2}}$$
(and $H^1(\cI_Z \otimes \cO(a)) = 0)$.
\smallskip

\noindent 
If we let $Z_m$ be the subspace defined by $\cI_Z^m,$ we obtain in
completely the same way (considering $\cI_Z^m \otimes \cO(ma)$ and
assuming $H^0(\cI_Z^m \otimes \cO(ma)) = 0$) that
$$
l(Z_m) = {{ma(ma+1)} \over {2}}.
$$
Now consider the case that $Z$ is reduced. Then $l(Z_2) = 3 l(Z) $
and $l(Z_3) = 10 l(Z).$ On the other hand, we can compute $l(Z_i)$ by
the above formula, and this produces a contradiction.  So $Z$ cannot
be reduced.\smallskip

\noindent
To deal with the general case, choose a deformation $(Z_t)$
of $Z = Z_0$ such that $Z_t$ is reduced for $t \ne 0.$ Then clearly
$\cI_{Z_t} \otimes \cO(a)$ is pseudo-effective. Now we apply to $Z_t$
the above considerations and obtain a contradiction.  \qed

\claim 9.7 Proposition| Let $X$ be a smooth projective 4-fold with
$K_X$ pseudo-effective. If $p(K_X) = 3$, then $\kappa(X) \geq 0$.
\endclaim

\proof. Here any reduction is an elliptic fibration.  We choose a
holomorphic birational model $f: X \la Y$ (with $X$ and $Y$ smooth),
such that {\item {(a)} $f$ is smooth over $Y_0$ and $Y \setminus Y_0$
is a divisor with simple normal crossings only;
\item{(b)} the $j-$function extends to a holomorphic map $J : Y \la \bP_1.$ }

By the first property, $f_*(K_X) $ is locally free [Ko86], and we
obtain the well-known formula of $\bQ-$divisors
$$
K_X = f^*(K_Y + \Delta) + E \leqno (*) $$
Here $E$ is an effective divisor such that
$f_*(\cO_X(E)) = \cO_Y$. Moreover
$$
\Delta = \Delta_1 + \Delta_2 $$
with $$\Delta_1 = \sum (1 - {{1} \over {m_i}}) F_i + \sum a_k
$$
and
$$
\Delta_2 \sim {{1} \over {12}} J^*(\cO(1)).
$$
Here $F_i$ are the components over which we have multiple fibers
and $D_k$ are the other divisor components over which there singular
fibers. The $a_k \in {{1} \over {12}} \bN$ according to Kodaira's
list.  Then by a general choice of the divisor $\Delta_2,$ the pair
$(Y,\Delta_1+\Delta_2)$ is klt. Now $K_Y + \Delta$ is
pseudo-effective. In fact, by Theorem 2.2 it suffices to show that
$K_Y + \Delta \cdot C_t \geq 0$ for every covering family of curves.
But this is checked very easily as in (9.4/9.6).  Hence the log
Minimal Model Program [Ko92] in dimension 3 implies that $K_Y +
\Delta$ is effective. Hence $\kappa (X) \geq 0$ by (*).  \qed

In order to attack the case $p(K_X) = 0$, i.e.\ there is a connecting
family $(C_t)$ with $K_X \cdot C_t = 0$, we prove a more general
result.

\claim 9.8 Theorem| Let $X$ be a projective manifold of any dimension
$n$, $(C_t)$ a connecting family and $L$ pseudo-effective. If $L \cdot
C_t = 0$, then there exists a line bundle $L'$ with $L \equiv L'$ and
$\kappa(L') = 0$. More generally, if $L \equiv N + Z$ is the
divisorial Zariski decomposition, then $Z = 0.$ \endclaim

\claim 9.9 Corollary| Let $X$ be a projective manifold with $q(X)=0$
such that $K_X$ is pseudo-effective. If $p(K_X) = 0$, i.e.\ there 
is a connecting family $(C_t)$ such that $K_X \cdot C_t = 0$, then 
$\kappa(X) = 0$.
\endclaim

First we derive 9.9 from 9.8.
\bigskip

\noindent
{\it Proof of 9.9.} Applying (9.8), $mK_X$ is effective up to some
numerically trivial line bundle $G$. Since we may assume $q(X) = 0$, 
the line bundle $G$ is trivial after finite \'etale cover,
hence $\kappa(X) \geq 0$. By 8.5 we get $\kappa(X) = 0$.  \qed
\bigskip

\noindent
{\it Proof of 9.8.}
\smallskip

\noindent First notice that suffices to prove the statement on the 
Zariski decomposition; then the
first statement follows from [Bou02b], see (3.7). In fact, then $L
\equiv N$ and $N$ is an effective $\bQ-$divisor.  \bigskip \noindent
(I) In a first step we assume that $\dim T = n-1$ resp.\ we can find a
{\it connecting} $(n-1)$-dimensional subfamily.  \smallskip

\noindent
Let 
$$
L = N_1 + Z_1
$$
be the divisorial Zariski decomposition of $L$ with $N_1$ effective 
and $Z_1$ nef in codimension 1 (as real divisors). \\
We consider the graph $p_0: \cC_0 \to X$ with projection $q_0: \cC_0
\to T;$ we may assume $T$ smooth and $\cC_0$ normal. Let
$\pi: \cC \to \cC_0$ be a desingularisation and put 
$p = p_0 \circ \pi$, $q = q_0 \circ \pi$. \\
Now $p^*(Z_1) $ might not be nef in codimension 1, so that we
decompose $p^*(Z_1) = \tilde N + Z$ and put $N = p^*(N_1) + \tilde N$.
We end up with
$$
p^*(L) = N + Z
$$
and this is the divisorial Zariski decomposition of $p^*(L)$. This
equation holds up to the pull-back of a topologically trivial
$\bR$-line bundle on $X$.  Notice that $N$ does not meet the general
fiber of $q$ since $L \cdot C_t = 0$.  This decomposition can be
rewritten as follows
$$
p^*(L)  =  \sum a_i F_i + \sum b_j B_j + q^*(N') + q^*(A) - E' \leqno (*)
$$
where $F_i$ are fibered over $q(F_i)$ by parts of reducible
1-dimensional fibers, the $B_j$ are irreducible components of the
exceptional locus of $\pi$ with $\codim q(B_j) \geq 2$, where $N'$ is
$\bR$-effective, where $E'$ is effective with ${\rm codim}(q(E')) \geq 2$
and finally $Z = q^*(A)$. In particular $A$ is nef in
codimension~$1$. The coefficients $a_i,b_j$ are positive and real a priori.
\bigskip

\noindent
First we verify the decomposition
$$
N  = \sum a_i F_i + \sum b_j B_j + q^*(N'). \leqno (**)
$$
In fact, since $p_0^*(N_1) \cdot q^{-1}(t) = 0$ for general $t$, we
can decompose the effective $\bR$ divisor $p_0^*N_1$ into the
components $F'_j $ with $ \codim q_0(F'_j) = 1$ but the $F'_j$ contain
only parts of fibers, and the other components which then consist only
of full fibers of $q_0$.  Therefore we can write
$$
p_0^*(N_1) = \sum a'_i F'_i + q_0^*(N')
$$
and then
$$
p^*N_1 = \sum a'_i F_i + \sum c_j B_j + q^*(N')$$
with $c_j \geq 0$.
Then $(**)$ follows by adding and decomposing $\tilde N$; notice here
that $\tilde N$ cannot contain multisections of $q$ since $p^*(L)
\cdot q^{-1}(t) = 0$.
\bigskip

\noindent  
To get the decomposition $(*)$, consider first a line bundle $B$ over
$\cC$ which is $q$-nef in codimension 1 with $B \cdot q^{-1}(t) = 0$
for general $t$. This means that $B$ is $q$-nef over a Zariski open
subset in $T$ whose complement has codimension at least 2.  Then $B$
defines a section in $R^1q_*(\cO_{\cC})$ over a Zariski open affine
set $T_0 \subset T$ such that ${\rm codim}(T \setminus T_0) \geq 2$.
Since $R^1q_*(\cO_{\cC})$ is a direct sum of a torsion sheaf supported
in codimension at least 2 and a reflexive sheaf ([Ko86]), $s$ extends
to a section of $R^1q_*(\cO_{\cC})$ on all of $T$ (possibly first enlarge
$T_0$).  By the Leray spectral sequence we obtain a topologically
trivial line bundle $G$ on $\cC$ such that $B \vert q^{-1}(t) = G
\vert q^{-1}(t) $ for general $t$. Thus, possibly substituting $B$ by
$B \otimes G^*$, we may assume that $B \vert q^{-1}(t) $ is trivial
for general $t$.  Now consider the canonical map $q^*q_*(B) \to B$ to
obtain a decomposition of type $(*)$ for $B$, namely $B = q^*(B') +
E$, where $E$ is a not necessarily effective divisor whose components
$E_k $ satisfy $\codim q(E_k) \geq 2$.
\smallskip

\noindent
Going back to our case, we would like to apply this to $B = Z$.
However $Z \in H^{1,1}_{\bQ}(\cC) \otimes \bR$ is not a $\bQ$-divisor.
Therefore we approximate $Z$ by $\bQ$-divisors $Z_j \in
H^{1,1}_{\bQ}(\cC)$ such that $Z_j \cdot q^{-1}(t) = 0$ for general
$t$. This is possible since the linear subspace $ \{ A \vert A \cdot
q^{-1}(t) = 0 \}$ is rationally defined. Now apply the previous
considerations to $Z_j;$ but first we have to make sure that the $Z_j$
are $q$-nef in codimension 1. This can be achieved by requiring that
the $Z_j$ are 0 on the components of general reducible fibers.  Now
applying our previous considerations, it follows that $Z_j = q^*(Z_j')
+ E_j$ with divisors $E_j$ whose components $E_{k,j} $ satisfy $\codim
q(E_{k,j}) \geq 2$.  Then let $Z'$ be the limit of the $Z_j;$ we obtain
$Z = q^*(Z') + E'$. Since $Z$ is nef in codimension 1, we conclude
that $-E'$ is effective.  This finally establishes $(*)$ (possibly we have to
pass from $L$ to $mL$ in order to avoid multiple components in $\sum a_i F_i$).
\bigskip

\noindent
We also notice that the coefficients $a_i$ in $(*)$ must be rational.
In fact, by $(*)$, $\sum a_i F_i$ is rational and by Boucksom [Bou02b,
2.1.15]; see (3.7), the $F_i$ are linearly independent.
\bigskip

\noindent
(I.a) First suppose that $\deg p = 1$. Then our claim comes down to prove that 
\bigskip

\noindent
$$ Z = 0. \eqno (*{*}*)$$  
\bigskip

\noindent

Since $\deg p = 1$, the map $p$ is birational; let $E$ denote the
exceptional locus.  Then two curves from the family $(C_t)$ can only
meet at points in $p(E)$.  Since $(C_t)$ is connecting, $E$ projects
onto $T$, i.e.\ $q(E) = T$. To be more precise, we pick some component
$E_i$ projecting onto $T$. Thus every $t$ is contained in some
subvariety
$$
T(x):= q(p^{-1}(x) \cap E_i),
$$
i.e.\ every $C_t$ passes through a
point $x \in p(E_i)$ and through each such $x$ there exists an at
least 1-dimensional subfamily $(C_t)$.  Notice that $T(x)$ might not
be irreducible, so that for general $x \in p(E_i)$ we pick an
irreducible component such that (taking closure) we obtain a compact
family which we again denote $(T(x))$.
\smallskip

\noindent
Now consider $x \in p(E_i)$ general and let $S_x \subset q^{-1}(T(x))$
be the irreducible component mapping onto $T(x)$.  By restricting
$$
p^*(L) = \sum a_iF_i + \sum b_j B_j + q^*(N'+A)
$$
to $S_x \subset q^{-1}(T(x))$ we obtain 
$$
A \vert T(x) \equiv N' \vert T(x) \equiv 0
$$ (cut by a general hyperplane section in $T(x)$).
If $T(x) = T$, then $N'$ and $A$ are numerically trivial and $(*{*}*)$
holds. 
\smallskip \noindent 
So suppose $\dim T(x) < \dim T$. If the family $(T(x))$ is
connecting, then we find connecting families of curves, say $(C_s)$
such that 
$$N'+A \cdot C_s = 0,$$ 
hence $N'+A $ is $\bQ$-effective by
induction and $(*{*}*)$ holds. 
\smallskip

\noindent
If $(T(x))$ is not connecting, then we form the quotient $g: T \merto
W$. We can choose $x,x' \in p(E)$ such that $T(x) \cap T(x') =
\emptyset.$ Since $C_{t_1} \cap C_{t_2} \subset p(E)$ for all choices
$t_j \in T$, we conclude that $x$ and $x'$ cannot be connected by
chains of $C_t$'s. This is a contradiction.

\bigskip
\noindent
(b) Now let $\deg p \geq 2$. Take a general $C_t$. Then through the
general $x \in C_t$ there is at least one other $C_{t'}$.  Therefore
we obtain a 1-dimensional family $(C_s)_{s \in T_t}$ through $C_t$. To
be more precise, let $D \subset X$ be the subspace over which $p$ has
positive dimensional fibers.  Then let
$$
T_t \subset q({\overline {p^{-1}(C_t \setminus D)}})
$$
be the union of all irreducible components of dimension 1 whose
$q-$images still have dimension 1.
Putting things together, we obtain a family $(T_t),$ however the
general $T_t$ might be reducible.  We claim that the family $(T_t)$ is
connecting.  In fact, take $t_1, t_2 \in T$ general. Then we can join
the curves $C_{t_1}$ and $C_{t_2}$ by an odd number of irreducible
curve $C_{s_j}, 1 \leq j \leq 2n-1.$ So $C_{t_1} \cap C_{s_1} \ne
\emptyset,$ $C_{s_1} \cap C_{s_2} \ne \emptyset$ etc. Thus $s_1 \in
T_{t_1} \cap T_{s_2}$ so that $T_{t_1} \cap T_{s_2} \ne \emptyset.$
Moreover $s_3 \in T_{s_2} \cap T_{s_4}$ so that $T_{s_2} \cap T_{s_4}
\ne \emptyset.$ Continuing, finally $T_{s_{2n-2}} \cap T_{t_2} \ne
\emptyset,$ so that $t_1$ and $t_2$ can be joined by chains of $T_s.$
Thus $(T_t)$ is connecting.  \smallskip \noindent Let $S_t \subset
\cC_t$ be the corresponding surface over $T_t$ resp. an irreducible
component.  Then $p^*(L) $ is numerically trivial on the general fiber
of $q \vert S_t$ and also on some multi-section. Therefore $q^*(N'+A)
\vert S_t = 0$ due to the following remark (+), hence $(N'+A) \cdot
T_t = 0$.  \smallskip

\noindent
Let $L = N'+A$ be a pseudo-effective line bundle over a smooth
projective surface or threefold $T$ with $N$ an effective
$\bR$-divisor and $A$ nef in codimension~$1$.  Assume that there is a
map $g: T \to W$ to the smooth curve $W$. Then $L$ is numerically
$\bQ$-effective.\smallskip

\noindent 
In fact, $N = \sum r_i F'_i$ with fiber components $F'_i$ and $A$ is
numerically trivial on all fibers since $A$ is nef in codimension~$1$.
Thus $A \equiv g^*(A')$ with $A'$ nef on $W$. This proves (+).
\smallskip \noindent
Now we apply induction if the general $T_t$ is irreducible resp. Lemma 9.10 if
the general $T_t$ is reducible and obtain that $A = 0, $ hence $Z = 0.$

\bigskip
\noindent
(II) If there is no connecting $(n-1)$-dimensional subfamily, we
choose some \hbox{$(n-1)$}-dimensional subfamily $\cC_0$ over $T_0 \subset
T$.  Now consider the $p$-preimages of the $C_t$ corresponding to $t
\not \in T_0$. Then we obtain a connecting family of multi-sections on
which $q^*(N'+A) $ is numerically trivial and, taking $q$-images, a
connecting family in $T_0$ on which $N'+A$ is numerically trivial.
Then argue by induction as in (I.b).\qed

\claim 9.10 Lemma| If Theorem 9.8 holds in dimension at most $n$, then
it also holds in dimension $\leq n$ for arbitrary covering connecting
families $(C_t)$ with the general $C_t$ being reducible: if $L$ is
pseudo-effective and $L \cdot C_t^j = 0$ on every movable component
$C_t,$ then $Z = 0$ in the divisorial Zariski decomposition.
\endclaim

\proof. Let $(C_t)_{t \in T}$ be a connecting family of reducible
curves with graph $p: \cC \to X$ and $q: \cC \to T.$ Let
$$
\cC = \bigcup_i^l \cC_i
$$
be the decomposition into irreducible components and set $p_i = p
\vert \cC_i, q_i = q \vert \cC_i.$ Some $\cC_i$ might still have
generically reducible $q$-fibers, so we pass to the normalisation
$\tilde \cC_i \la \cC_i.$ Let $\tilde \cC_{i,j}$ be the decomposition
into irreducible components; then every $\tilde \cC_{i,j}$ defines a
family of generically irreducible curves on $X.$ In total we obtain
finitely many families $\hat \cC_ k \la S_k$ of generically
irreducible curves, not all of them covering possibly.  \smallskip
\noindent Pick one of the covering families, say $\hat \cC_1.$ Let $f:
X \merto X_1$ be the associated quotient (with $X_1$ smooth). Since we
are allowed to blow up, we may assume $f$ holomorphic from the
beginning.  If $\hat \cC_1$ happens to be connecting, i.e. $\dim X_1 =
0,$ then (9.8) gives our claim for the given line bundle $L.$ So
suppose $\dim X_1 > 0.$ Our plan is to proceed by induction on $\dim
X.$ Thus consider the induced family $\tilde C_t = f_*(C_t)$ (for $t$
generic; then take closure in the cycle space).  Obviously $(\tilde
C_t)$ is connecting.  By (9.8) and (8.5) - possibly after tensoring
with a topologically trivial line bundle as in the proof of (9.8) - we
have $\kappa (L \vert F) = 0$ for suitable $m$ and the general fiber
$F$ of $f.$ Thus $f_*(mL)$ has rank $1.$ Choosing $m$ sufficiently
divisible, we obtain
$$
mL = f^*(L_1) + E_1 - E_2,
$$
where $E_j$ are effective, the components of $E_1$ either consists
of parts of fibers of $f$ or project onto subvarieties of codimension
at least 2 in $X_1$ and the components of $E_2$ also project onto
subvarieties of codimension at 2 in $X_1.$ By passing to a suitable
model of $f$ (blow up $X$ and $Y$), we may assume that $f_*(mL) $ is
locally free, hence $E_2 = 0.$ Now by cutting down to movable curves,
we conclude as in 9.4/9.6 that $L_1$ is pseudo-effective.  In order to
apply induction we still need to show that $L \cdot \tilde C_t^j = 0$
for the movable components.  This is however clear:
$$ 0 = mL \cdot C_t^j = E_1 \cdot C_t^j + L' \cdot \tilde C_t^j, $$
so that $E_1 \cdot C_t^j = L' \cdot \tilde C_t^j = 0.$ 
\smallskip

\noindent Hence we can apply induction: if $L_1 = N_1 + Z_1$ is the
divisorial Zariski decomposition, then
$Z_1 = 0.$ Thus the same holds for $L.$ \qed

Combining everything in this section we finally obtain

\claim 9.11 Theorem| Let $X$ be a smooth projective 4-fold (or a normal
projective 4-fold with only canonical singularities).  If $K_X$ is
pseudo-effective and if there is a covering family $(C_t)$ of curves
such that $K_X \cdot C_t = 0$, then $\kappa(X) \geq 0$.
\endclaim

The remaining task is to consider 4-folds $X$ with $K_X$
pseudo-effective such that $K_X \cdot C > 0$ for every curve $C$
passing through a very general point of $X$, i.e.\ $p(K_X) = 0$.  In
that case one expects that $X$ is of general type. It is easy to see
that every proper subvariety $S$ of $X$ passing through a very general
point of $X$ is of general type, i.e.\ its desingularisation is of
general type; see 9.12 below. But it is not at all clear whether $K_X
\vert S$ is big, which is of course still not enough to conclude. 

\claim 9.12 Proposition| Let $X$ be a smooth projective 4-fold with
$p(K_X) = 0$. Then every proper subvariety $S \subset X$ through a
very general point of $X$ is of general type.
\endclaim

\proof.  Supposing the contrary, we find a covering family $(S_t)$ of
subvarieties such that the general $S_t$, hence every $S_t$, is not of
general type. Consider the desingularised graph $p: \cC \to X$ of
this family; by passing to a subfamily we may assume $p$ generically
finite.  Denoting $q: \cC \to T$ the parametrising projection, the
general fiber $\hat S_t$ is a smooth variety of dimension at most 3
and not of general type.  Using a minimal model if $\kappa(\hat S_t)
= 0$, we find in $\hat S_t$ a covering family of curves intersection
$K_{\hat S_t} $ trivially.  Thus we find a covering family $(C_s)$ in
$\cC$, all members being in $q$-fibers such that $K_{\cC} \cdot C_s =
0$. Since $K_{\cC} = p^ *(K_X) + E$ with $E$ effective, we get $K_X
\cdot p_*(C_s) \leq 0$, a contradiction.
\qed

Using the Iitaka fibration we obtain

\claim 9.13 Proposition| Let $X$ be a smooth projective 4-fold with
$p(K_X) = 0$. Then $\kappa(X) \ne 1,2,3$.
\endclaim

\section{\S10 Appendix: towards transcendental Morse inequalities}

As already pointed out, for the general case of the conjecture 2.3
a transcendental version of the holomorphic Morse inequalities 
would be needed. The expected statements are contained in the following
conjecture

\claim 10.1 Conjecture| Let $X$ be a compact complex manifold,
and $n=\dim X$.
\smallskip
\item{\rm(i)} Let $\alpha$ be a closed, $(1, 1)$-form on $X$. We denote by
$X(\alpha, \leq 1)$ the set of points $x\in X$ such that $\alpha_x$
has at most one negative eigenvalue. If 
$\int_{X(\alpha, \leq 1)}\alpha^n>0$, the class $\{\alpha\}$ contains
a K\"ahler current and
$$\Vol(\alpha)\geq \int_{X(\alpha, \leq 1)}\alpha^n.$$
\smallskip
\item{\rm(ii)} Let $\{\alpha\}$ and $\{\beta\}$ be nef cohomology classes
of type $(1,1)$ on $X$ satisfying the inequality 
$\alpha^n-n\alpha^{n-1}\cdot\beta>0$. Then $\{\alpha-\beta\}$ contains 
a K\"ahler current and
$$\Vol(\alpha-\beta)\geq \alpha^n - n\alpha^{n-1}\cdot\beta.$$
\vskip0pt
\endclaim

\noindent{\bf Remarks about the conjecture.} If $\alpha= c_1(L)$ for some 
holomorphic line bundle $L$ on
$X$, then the inequality $(**)$ was established in [Bou02a] as a consequence
of the results of [De85]. In general, (ii) is a consequence of (i). In fact,
if $\alpha$ and $\beta$ are smooth positive definite $(1,1)$-forms and
$$\lambda_1\geq\ldots\geq\lambda_n>0$$
are the eigenvalues of $\beta$ with respect to $\alpha$, then
$X(\alpha-\beta,\leq 1)=\{x\in X\,;\;\lambda_2(x)<1\}$ and
$${\bf 1}_{X(\alpha-\beta,\leq 1)}(\alpha-\beta)^n=
{\bf 1}_{X(\alpha-\beta,\leq 1)}(1-\lambda_1)\ldots(1-\lambda_n)\ge 
1-(\lambda_1+\ldots+\lambda_n)$$
everywhere on $X$. This is proved by an easy induction on $n$. An integration
on $X$ yields inequality (ii). In case $\alpha$ and $\beta$ are just nef but
not necessarily positive definite, we argue by considering
$(\alpha+\varepsilon\omega)-(\beta+\varepsilon\omega)$ with a positive
hermitian form $\omega$ and $\varepsilon>0$ small.\qed
\medskip

The full force of the conjecture is not needed here. First of all,
we need only the case when $X$ is compact K\"ahler. Let us consider a big 
class $\{\alpha\}$, and a sequence of K\"ahler currents
$T_m\in \{\alpha\}$ with logarithmic poles, such that
there exists a modification $\mu_m:X_m\mapsto X$, with the properties

{\itemindent=1.3cm
\item {$(10.2')$} $\mu_m^*T_m= \beta_m+ [E_m]$ where $\beta_m$ is a
  semi-positive $(1, 1)$-form, and $E_m$ is an effective $\bQ$-divisor
  on $X_m$.  \smallskip
\item {$(10.2'')$} 
$\Vol (\{\alpha\})= \lim_{m\mapsto \infty}\int_X\beta_m^n$.
\smallskip
\noindent
(see Definition 3.2).}

A first trivial observation is that the following 
uniform upper bound for $c_1(E_m)$ holds.

\medskip
\claim 10.3 Lemma| Let $\omega$ be a 
K\"ahler metric on $X$, such that $\{\omega- \alpha\}$ contains a smooth, 
positive representative. Then for each $m\in \bZ_+$, the $(1, 1)$-class 
$\mu_m^*\{\omega\}- c_1(E_m)$ on $X_m$ is nef.
\endclaim

\medskip

\proof. If $\gamma$ is a smooth positive representative in 
$\{\omega-\alpha\}$, then $\mu_m^\star\gamma+\beta_m$ is a smooth
semi-positive representative of $\mu_m^*\{\omega\}- c_1(E_m)$.\qed

\bigskip A second remark is that in order to prove the duality
statement 2.3 for {\sl projective} manifolds, it is enough to
establish the estimate
$$\Vol (\omega-  A)\geq \int_X\omega^n- n\int_X\omega^{n- 1}
\wedge c_1(A)\leqno (*)$$
where $\omega$ is a K\"ahler metric, and $A$ is an ample line bundle on $X$.
Indeed, if $\{\alpha\}$ is a big cohomology class, we use the above
notations and we can write
$$\beta_m+ tE_m= \beta_m+ t\mu_m^*A- t(\mu_m^*A- E_m)$$
where $A$ is an ample line bundle on $X$ such that $c_1(A)- \{\alpha\}$ 
contains a smooth, positive representative. The arguments of the proof of 4.1
will give the orthogonality estimate, provided that we are able to
establish $(*)$.

\bigskip 
In this direction, we can get only a weaker statement with a
suboptimal constant~$c_n$.

\medskip
 
\claim 10.4 Theorem (analogue of Lemma 4.2)| Let $X$ be a projective
manifold of dimension $n$.  Then there exists a constant $c_n$
depending only on dimension $($actually one can take
$c_n=(n+1)^2/4\,)$, such that the inequality
$$\Vol (\omega-A)\geq \int_X\omega^n- c_n\int_X\omega^{n- 1}\wedge c_1(A)$$
holds for every K\"ahler metric $\omega$ and every ample line bundle
$A$ on~$X$.
\endclaim

\proof. Without loss of generality, we can assume that $A$ is very ample
(otherwise multiply $\omega$ and $A$ by a large positive integer).
Pick generic sections $\sigma_0,\sigma_1,\ldots,\sigma_n\in|A|$ so that
one gets a finite map
$$
F:X\to\bP_\bC^n,\qquad x\mapsto[\sigma_0(x):\sigma_1(x):\ldots:\sigma_n(x)].
$$
We let $\theta=F^*\omega_\FS\in c_1(A)$ be the pull-back of the Fubini-Study
metric on $\bP_\bC^n$ (in particular $\theta\ge 0$ everywhere on $X$), and put
$$
\psi=\log{|\sigma_0|^2\over|\sigma_0|^2+|\sigma_1|^2+\ldots+|\sigma_n|^2}.
$$
We also use the standard notation $d^c={i\over4\pi}(\ol\partial-\partial)$
so that $dd^c={i\over 2\pi}\ddbar$. Then 
$$
dd^c\psi=[H]-\theta
$$ 
where $H$ is the hyperplane section $\sigma_0=0$ and $[H]$ is the current
of integration over $H$ (for simplicity, we may further assume that $H$ is
smooth and reduced, although this is not required in what follows). The set 
$U_\varepsilon=
\{\psi\le 2\log\varepsilon\}$ is an $\varepsilon$-tubular neighborhood of $H$.
Take a convex increasing function $\chi:\bR\to\bR$ such that $\chi(t)=t$
for $t\ge 0$ and $\chi(t)={}$constant on some interval $]-\infty,t_0]$.
We put $\psi_\varepsilon=\psi-2\log\varepsilon$ and
$$
\alpha_\varepsilon:=dd^c\chi(\psi_\varepsilon)+\theta=
(1-\chi'(\psi_\varepsilon))\theta+\chi''(\psi_\varepsilon)d\psi_\varepsilon
\wedge d^c\psi_\varepsilon\ge 0.
$$
Thanks to our choice of $\chi$, this is a smooth form with support
in $U_\varepsilon$. In particular, we find
$$
\int_{U_\varepsilon}\alpha_\varepsilon^n=
\int_{U_\varepsilon}\alpha_\varepsilon\wedge\theta^{n-1}=
\int_X\theta^n=c_1(A)^n.
$$
It follows from these equalities that we have $\lim_{\varepsilon\to 0}
\alpha_\varepsilon=[H]$ in the weak topology of currents.
Now, for each choice of positive parameters 
$\varepsilon, \delta$, we consider the Monge-Amp\`ere equation
$$
(\omega+ i\ddbar \varphi_\varepsilon)^n
=(1-\delta)\omega^n+\delta {{\int_X \omega^n}\over
c_1(A)^n}\;\alpha_\varepsilon^n.\leqno(10.5)
$$
By the theorem of S.-T.~Yau [Yau78], there exists a smooth solution
$\varphi_\varepsilon$, unique up to normalization by an additive
constant, such that $\omega_\varepsilon:= \omega+ i\ddbar
\varphi_\varepsilon> 0$. Since
$\int_X\omega_\varepsilon\wedge\omega^{n-1}=\int_X\omega^n$ remains
bounded, we can extract a weak limit $T$ out of the family
$\omega_\varepsilon\,$; then $T$ is a closed positive current, and the
arguments in [Bou02a] show that its absolutely continuous part
satisfies
$$\int_XT_{ac}^n\geq (1- \delta)\int_X\omega^n.$$
We are going to use the same ideas as in [DPa03], in order to estimate the
singularity of the current $T$ on the hypersurface $H$. For this,
we estimate the integral $\int_{U_\varepsilon}\omega_\varepsilon \wedge
\theta^{n-1}$ on the tubular neighborhood $U_\varepsilon$ of $H$.
Let us denote by $\rho_1\le\ldots\le\rho_n$ the eigenvalues of 
$\omega_\varepsilon$ with respect to $\alpha_\varepsilon$, computed
on the open set $U'_\varepsilon\subset U_\varepsilon$ where
$\alpha_\varepsilon$ is positive definite. The Monge-Amp\`ere equation
(10.5) implies
$$
\rho_1\rho_2\ldots\rho_n\ge \delta{\int_X\omega^n\over c_1(A)^n}.
$$
On the other hand, we find $\omega_\varepsilon\ge\rho_1\alpha_\varepsilon$
on $U'_\varepsilon$, hence
$$
\int_{U_\varepsilon}\omega_\varepsilon\wedge\theta^{n-1}
\ge \int_{U'_\varepsilon}\rho_1\alpha_\varepsilon\wedge\theta^{n-1}
\ge \delta{\int_X\omega^n\over c_1(A)^n}
\int_{U'_\varepsilon}{1\over\rho_2\ldots\rho_n}
\alpha_\varepsilon\wedge\theta^{n-1}.
\leqno(10.6)
$$
In order to estimate the last integral in the right hand side, we
apply the Cauchy-Schwarz inequality to get
$$
\Big(\int_{U'_\varepsilon}(\alpha_\varepsilon^n)^{1/2}
(\alpha_\varepsilon\wedge\theta^{n-1})^{1/2}\Big)^2\le
\int_{U'_\varepsilon}\rho_2\ldots\rho_n\alpha_\varepsilon^n
\int_{U'_\varepsilon}{1\over\rho_2\ldots\rho_n}\alpha_\varepsilon\wedge
\theta^{n-1}.\leqno(10.7)
$$
By definition of the eigenvalues $\rho_j$, we have
$$
\int_{U'_\varepsilon}\rho_2\ldots\rho_n\alpha_\varepsilon^n\le
n\int_X\omega_\varepsilon^{n-1}\wedge\alpha_\varepsilon=
n\int_X\omega^{n-1}\wedge c_1(A).\leqno(10.8)
$$
On the other hand, an explicit calculation shows that
$$
\eqalign{
&\alpha_\varepsilon^n\ge n(1-\chi'(\psi_\varepsilon))^{n-1}
\chi''(\psi_\varepsilon)\,d\psi_\varepsilon\wedge d^c\psi_\varepsilon
\wedge\theta^{n-1},\cr
&\alpha_\varepsilon\wedge\theta^{n-1}\ge
\chi''(\psi_\varepsilon)\,d\psi_\varepsilon\wedge d^c\psi_\varepsilon
\wedge\theta^{n-1},\cr}
$$
hence
$$
\int_{U'_\varepsilon}(\alpha_\varepsilon^n)^{1/2}
(\alpha_\varepsilon\wedge\theta^{n-1})^{1/2}\ge
n^{1/2}\int_X
(1-\chi'(\psi_\varepsilon))^{(n-1)/2}\chi''(\psi_\varepsilon)
\,d\psi_\varepsilon\wedge d^c\psi_\varepsilon
\wedge\theta^{n-1}
$$
(we can integrate on $X$ since the integrand is zero anyway outside
$U'_\varepsilon$). Now, we have
$$
\eqalign{
{n+1\over 2}&(1-\chi'(\psi_\varepsilon))^{(n-1)/2}\chi''(\psi_\varepsilon)
\,d\psi_\varepsilon\wedge d^c\psi_\varepsilon\cr
&=-d\Big((1-\chi'(\psi_\varepsilon))^{(n+1)/2}d^c\psi_\varepsilon\Big)+
(1-\chi'(\psi_\varepsilon))^{(n+1)/2}dd^c\psi_\varepsilon\cr
&=-d\Big((1-\chi'(\psi_\varepsilon))^{(n+1)/2}d^c\psi_\varepsilon\Big)+
[H]-(1-\chi'(\psi_\varepsilon))^{(n+1)/2}\theta\cr}
$$
and from this we infer
$$
\eqalign{
{n+1\over 2}&\int_X(1-\chi'(\psi_\varepsilon))^{(n-1)/2}
\chi''(\psi_\varepsilon)
\,d\psi_\varepsilon\wedge d^c\psi_\varepsilon\wedge\theta^{n-1}\cr
&=\int_X[H]\wedge\theta^{n-1}-\int_X
(1-\chi'(\psi_\varepsilon))^{(n+1)/2}\theta^n\cr
&\to c_1(A)^n\qquad\hbox{as $\varepsilon\to 0$}.\cr}
$$
We thus obtain
$$
\int_{U'_\varepsilon}(\alpha_\varepsilon^n)^{1/2}
(\alpha_\varepsilon\wedge\theta^{n-1})^{1/2}\ge{2\sqrt{n}\over n+1}c_1(A)^n
-o(1)\qquad\hbox{as $\varepsilon\to 0$}.\leqno(10.9)
$$
The reader will notice, and this looks at first a bit surprising, that
the final lower bound does not depend at all on the choice of $\chi$.
This seems to indicate that our estimates are essentially optimal
and will be hard to improve. Putting together (10.7), (10.8) and (10.9) 
we find the lower bound
$$
\int_{U'_\varepsilon}{1\over\rho_2\ldots\rho_n}\alpha_\varepsilon\wedge
\theta^{n-1}\ge
{4\delta\over(n+1)^2}{(c_1(A)^n)^2\over\int_X\omega^{n-1}\wedge c_1(A)}-o(1).
\leqno(10.10)
$$
Finally, (10.6) and (10.10) yield
$$
\int_{U_\varepsilon}\omega_\varepsilon\wedge\theta^{n-1}\ge
{4\delta\over(n+1)^2}{\int_X\omega^n\over\int_X\omega^{n-1}\wedge c_1(A)}
c_1(A)^n-o(1).
$$
As $\bigcap U_\varepsilon=H$, the standard support theorems for currents 
imply that the weak limit $T=\lim\omega_\varepsilon$ carries a divisorial 
component $c[H]$ with
$$
\int_X c[H]\wedge\theta^{n-1}\ge
{4\delta\over(n+1)^2}{\int_X\omega^n\over\int_X\omega^{n-1}
\wedge c_1(A)}c_1(A)^n.
$$
Therefore, as $[H]\equiv\theta\in c_1(A)$, we infer
$$
c\ge{4\delta\over(n+1)^2}{\int_X\omega^n\over\int_X\omega^{n-1}\wedge c_1(A)}.
$$
The difference $T-c[H]$ is still a positive current and has the same absolutely
continuous part as $T$. Hence
$$
\Vol(T-c[H])\ge\int_XT_{ac}^n\ge(1-\delta)\int_X\omega^n.
$$
The specific choice 
$$
\delta={(n+1)^2\over 4}{\int_X\omega^{n-1}\wedge c_1(A)\over\int_X\omega^n}
$$
gives $c\ge 1$, hence
$$
\Vol(T-[H])\ge\int_X\omega^n-{(n+1)^2\over 4}\int_X\omega^{n-1}\wedge c_1(A).
$$
Theorem 10.4 follows from this estimate.\qed

\claim 10.11 Remark| {\rm By using similar methods, we could also obtain an
estimate for the volume of the difference of two K\"ahler classes on
a general compact K\"ahler manifold, by using the technique of
concentrating the mass on the diagonal of $X\times X$ (see [DPa03]).
However, the constant $c$ implied by this technique also depends on
the curvature of the tangent bundle of $X$.}
\endclaim

We show below that the answer to conjecture 10.1 is positive at least 
when $X$ is a compact
hyperk\"ahler manifold (${}={}$compact irreducible holomorphic 
symplectic manifold). The same proof would work for a compact K\"ahler
manifold which is a limit by deformation of projective manifolds with 
Picard number $\rho=h^{1,1}$.

\claim 10.12 Theorem| Let $X$ be a compact hyperk\"ahler manifold,
and let $\alpha$ be a closed, $(1, 1)$-form on $X$. Then we have
$$\Vol(\alpha)\geq \int_{X(\alpha, \leq 1)}\alpha^n.$$
\endclaim

\proof. We follow closely the approach of D.~Huybrechts in
[Huy02], page 44. Consider $\cX\mapsto \Def(X)$ the universal
deformation of $X$, such that $\cX_0= X$. If $\beta\in H^2(X, \bR)$ is
a real cohomology class, then we denote by $S_\beta$ the set of points
$t\in \Def(X)$ such that the restriction $\displaystyle \beta_{\vert
  \cX_t}$ is of $(1, 1)$-type.

Next, we take a sequence of rational classes $\{\alpha_k\}\in H^2(X,
\bQ)$, such that $\alpha_k\to \alpha$ on $\cX$ as $k\mapsto
\infty$.  As $\{\alpha_k\}\to\{\alpha\}$, the hypersurface
$\displaystyle S_{\alpha_k}$ converge to $S_\alpha$; in particular, we
can take $\displaystyle t_k\in S_{\alpha_k}$ such that $t_k\to 0$.
In this way, the rational $(1, 1)$-forms $\displaystyle \alpha_{k
\vert \cX_{t_k}}$ will converge to our form $\alpha$ on $X$.

We have 
$$\eqalign{ \Vol(\alpha)\geq & \lim\sup_{k\mapsto
    \infty}\Vol(\alpha_{k \vert \cX_{t_k}})\geq \cr \geq &
  \lim\sup_{k\mapsto \infty}\int_{\cX_{t_k}(\alpha_{t_k}, \leq
    1)}\alpha_{t_k}^n \cr = & \int_{X(\alpha, \leq 1)}\alpha^n\cr }$$
where the first inequality is a consequence of the semi-continuity of
the volume obtained in [Bou02b], and the second one is a consequence of
the convergence statement above.

\claim 10.13 Corollary| If $X$ be a compact hyperk\"ahler manifold, or more
generally, a limit by deformation of projective manifolds with Picard
number $\rho=h^{1,1}$, then the cones $\cE$ and $\cM$ are dual.
\endclaim

\section{References}

\bigskip

{\eightpoint
  
\bibitem [Bou02a]&Boucksom, S.:& On the volume of a big line bundle;&
  Intern.\ J.\ Math.\ {\bf 13} (2002), 1043--1063&

\bibitem [Bou02b]&Boucksom, S.:& C\^ones positifs des vari\'et\'es
  complexes compactes;& Thesis, Grenoble 2002&

\bibitem[Ca81]&Campana, F.:& Cor\'eduction alg\'ebrique d'un espace
  analytique faiblement k\"ahl\'erien compact;& Inv.\ Math.\ {\bf 63}
  (1981) 187--223&
  
\bibitem[Ca94]&Campana, F.:& Remarques sur le rev\^etement universel des
  vari\'et\'es K\"ahleriennes compactes.& Bull.\ Soc.\ Math.\ 
  France.\ {\bf 122} (1994) 255-284&

\bibitem[CP91]&Campana, F.; Peternell, Th.:& Algebraicity of the ample
  cone of projective varieties.& J.\ f.\ reine u.\ angew.\ Math.\ {\bf
  407} (1991) 160-166&

\bibitem [CP01]&Campana, F.; Peternell, Th.:& The Kodaira dimension of
  Kummer threefolds.& Bull.\ Soc.\ Math.\ France {\bf 129}, 357-359
  (2001)&
 
\bibitem[CT86]&Colliot-Th\'el\`ene, J.L.:& Arithm\'etique des
  vari\'et\'es rationnelles et probl\`emes birationnels;& Proc.\ Intl.\ 
  Cong.\ Math.\ 1986, 641-653&
  
\bibitem [De85]&Demailly, J.-P.:& Champs magn\'etiques et in\'egalit\'es
  de Morse pour la $d''$-coho\-mo\-logie;& Ann.\ Inst.\ Fourier
  (Grenoble) {\bf 35} (1985) 189--229&

\bibitem [De90]&Demailly, J.-P.:& Singular hermitian metrics on
  positive line bundles;& Proceedings of the Bayreuth conference
  ``Complex algebraic varieties'', April~2-6, 1990, edited by
  K.~Hulek, T.~Peternell, M.~Schneider, F.~Schreyer, Lecture Notes in
  Math.\ ${\rm n}^\circ\,$1507, Springer-Verlag, 1992, 87--104&

\bibitem [De92]&Demailly, J.-P.:& Regularization of closed
  positive currents and intersection theory;& J.\ Alg.\ Geom.\ {\bf 1} (1992),
  361-409&
  
\bibitem [DEL00]&Demailly, J.P.; Ein, L.; Lazarsfeld,R.:& A
  subadditivity property of multiplier ideals.& Michigan Math.\ J.\ 
  {\bf 48} (2000), 137-156 &

\bibitem [DPa03]&Demailly, J.P.; Paun, M.:& Numerical characterization of
  the K\"ahler cone of a compact K\"ahler manifold.& math.AG/0105176 (2001),
  to appear in Annals of Math., 2003 &

\bibitem [DPS94]&Demailly, J.-P.; Peternell, Th.; Schneider, M.:&
  Compact complex manifolds with numerically effective tangent
  bundles;& J.\ Alg.\ Geom.\ {\bf 3} (1994) 295--345&

\bibitem [DPS96]&Demailly, J.-P.; Peternell, Th.; Schneider, M.:&
  Holomorphic line bundles with partially vanishing cohomology;&
  Israel Math.\ Conf.\ Proc.\ vol.~{\bf 9} (1996) 165--198&

\bibitem [DPS00]&Demailly, J.-P.: Peternell, Th.; Schneider, M.;& 
  Pseudo-effective line bundles on compact K\"ahler manifolds;& 
  Intern.\ J.\ Math. {\bf 6} (2001) 689-741&
  
\bibitem [Ec02]&Eckl, T.:& Tsuji's numerical trivial fibrations.& 
  math.AG/0202279 (2002)&
  
\bibitem [Fuj94]&Fujita, T.;& Approximating Zariski decomposition of
  big line bundles;& Kodai Math.\ J.\ {\bf 17} (1994), 1--3&

\bibitem [GHS03]&Graber, T.; Harris, J.; Starr, J.:& Families of
  rationally connected varieties;& Journal AMS {\bf 16} (2003), 57-67&

\bibitem [KMM87]&Kawamata, Y.; Matsuki, K.; Matsuda, K.:& Introduction
  to the minimal model program;& Adv.\ Stud.\ Pure Math.\ {\bf 10}
  (1987), 283-360&

\bibitem [KM92]&Koll\'ar, J.; Mori, S.:& Classification of
  three-dimensional flips;& Journal AMS {\bf 5} (1992), 533-703&

\bibitem [Ko86]&Koll\'ar, J.:& Higher direct images of dualizing
  sheaves II ;& Ann. Math. {\bf 124}, (1986), 171-202&

\bibitem [Ko87]&Kobayashi, S.:& Differential geometry of complex vector 
  bundles;& Princeton Univ. Press, 1987&
  
\bibitem[Ko92]&Kolla\'r, J.\ et al.:& Flips and abundance for algebraic
  $3$-folds;& Ast\'erisque {\bf 211}, Soc. Math. France 1992&
  
\bibitem [La00]&Lazarsfeld, R.:& Ampleness modulo the branch locus of
  the bundle associated to a bran\-ched covering;& Comm. Alg.{\bf 28}
  (2000), 5598-5599&
  
\bibitem [Mi87]&Miyaoka, Y.:& The Chern classes and Kodaira dimension
  of a minimal variety;& Adv. Stud. Pure Math. {\bf 10}, 449 - 476 (1987) &

\bibitem [Mi88]&Miyaoka, Y.:& On the Kodaira dimension of minimal
  threefolds;& Math. Ann. {\bf 281}, 325-332 (1988)&
  
\bibitem [Mi88a]&Miyaoka, Y.:& Abundance conjecture for 3-folds: case
  $\nu = 1$;& Comp.\ Math. {\bf 68}, 203-220 (1988)&

\bibitem [MM86]&Miyaoka, Y.; Mori, S.:& A numerical criteria for
  uniruledness;&  Ann.\ Math.\ {\bf 124} (1986) 65--69&
  
\bibitem [Mo87]&Mori, S.:& Classification of higherdimensional
  varieties;& Proc.\ Symp.\ Pure Math.\ {\bf 46} 269-331 (1987)&

\bibitem [Mo88]&Mori, S.:& Flip theorem and existence of minimal
  models of $3$-folds;& J.\ Amer.\ Math.\ Soc.\ {\bf 1} (1988) 177--353&

\bibitem [Na87]&Nakayama, N.:&On Weierstrass models;& Alg.\ Geom.\ and
  Comm.\ Algebra.\ Vol.\ in honour of Nagata, vol.~2, 
  Kinokuniya, Tokyo 1988, 405--431&

\bibitem [Pa98]&Paun, M.:& Sur l'effectivit\'e num\'erique des
  images inverses de fibr\'es en droites;& Math.\ Ann.\ {\bf 310}
  (1998) 411--421&
  
\bibitem [Pe01]&Peternell, Th.:& Towards a Mori theory on compact
  K\"ahler $3$-folds, III;& Bull.\ Soc.\
  Math.\ France {\bf 129}, 339-356 (2001)&

\bibitem [PSS99]&Peternell, Th.; Schneider, M.; Sommese, A.J.:& Kodaira 
  dimension of subvarieties.& Intl.\ J.\ Math.\ {\bf 10} (1999) 1065--1079&
  
\bibitem [PS02]&Peternell, Th.; Sommese, A.J.:& Ample vector bundles
  and branched coverings II.& To appear in the volume of the Fano
  conference, Torino 2002&

\bibitem [SB92]&Shepherd-Barron, N.:& Miyaoka's theorems on the generic
  semi-negativity of $T_X$.& Ast\'e\-risque {\bf 211} (1992) 103--114&

\bibitem [Ts00]&Tsuji, H.:&Numerically trivial fibrations.& math.AG/0001023 
  (2000)&

\bibitem [Workshop]&Bauer, T.\ et al.:& A reduction map for nef line 
  bundles;& Complex Geometry, vol.\ in honour of H.~Grauert, 
  Springer (2002), 27--36&

\bibitem [Yau78]&Yau, S.-T.:& On the Ricci curvature of a 
  complex K\"ahler manifold and the complex Monge--Amp\`ere equation;&
  Comm.\ Pure Appl.\ Math.\ {\bf 31} (1978), 339--411& 

}

\bigskip
\noindent
(version of May 14, 2004, printed on \today)
\bigskip\bigskip
{\parindent=0cm
S\'ebastien Boucksom, {\it Universit\'e de Paris VII}, 
boucksom@math.jussieu.fr\\
Jean-Pierre Demailly, {\it Universit\'e de Grenoble I}, 
demailly@fourier.ujf-grenoble.fr\\
Mihai Paun, {\it Universit\'e de Strasbourg I}, 
paun@math.u-strasbg.fr\\
Thomas Peternell, {\it Universit\"at Bayreuth}, 
Thomas.Peternell@uni-bayreuth.de

}
}
\end